\documentclass[reqno,11pt]{amsart} 
\usepackage{amsmath}
\usepackage{amssymb}
\usepackage{amsthm}

\newcommand\NoBlackBoxes{\global\overfullrule0pt}
\NoBlackBoxes
\theoremstyle{plain} 

\def\4{\kern1pt}

\def\6{\vphantom0}

\def\8{\kern-10pt}
\def\7#1{_{(#1)}}

\begin{document}

\title{POINCAR\'E INEQUALITIES AND NORMAL \\ 
APPROXIMATION FOR WEIGHTED SUMS}

\author{S. G. Bobkov$^{1,4}$}
\thanks{1) 
School of Mathematics, University of Minnesota, USA}

\author{G. P. Chistyakov$^{2,4}$}
\thanks{2) 
Faculty of Mathematics, University of Bielefeld, Germany}

\author{F. G\"otze$^{3,4}$}
\thanks{3) 
Faculty of Mathematics, University of Bielefeld, Germany}
\thanks{4) Research was supported by SFB 1283, Simons Foundation, 
and NSF grant DMS-1855575}


\begin{abstract}
Under Poincar\'e-type conditions, upper bounds are explored for 
the Kolmogorov distance between the distributions of weighted sums 
of dependent summands and the normal law. Based on improved 
concentration inequalities on high-dimensional Euclidean spheres, 
the results extend and refine previous results to
non-symmetric models.
\end{abstract}

\maketitle
\markboth{S. G. Bobkov, G. P. Chistyakov and F. G\"otze}{
CLT for weighted sums}

\def\theequation{\thesection.\arabic{equation}}
\def\E{{\mathbb E}}
\def\R{{\mathbb R}}
\def\C{{\mathbb C}}
\def\P{{\mathbb P}}
\def\I{{\mathbb I}}
\def\S{{\mathbb S}}

\def\s{{\mathfrak s}}

\def\H{{\rm H}}
\def\Im{{\rm Im}}
\def\Tr{{\rm Tr}}

\def\k{{\kappa}}
\def\M{{\cal M}}
\def\Var{{\rm Var}}
\def\Ent{{\rm Ent}}
\def\O{{\rm Osc}_\mu}

\def\ep{\varepsilon}
\def\phi{\varphi}
\def\F{{\cal F}}
\def\L{{\cal L}}

\def\be{\begin{equation}}
\def\en{\end{equation}}
\def\bee{\begin{eqnarray*}}
\def\ene{\end{eqnarray*}}


\section{{\bf Introduction}}
\setcounter{equation}{0}

\noindent
Let $X = (X_1,\dots,X_n)$ be an isotropic random vector in $\R^n$
($n \geq 2$), meaning that $\E X_i X_j = \delta_{ij}$ for all $i,j \leq n$,
where $\delta_{ij}$ is the Kronecker symbol. Define the weighted sums
$$
S_\theta = \theta_1 X_1 + \dots + \theta_n X_n, \qquad 
\theta = (\theta_1,\dots,\theta_n), \ \ \theta_1^2 + \dots + \theta_n^2=1,
$$
with coefficients from the unit sphere $\S^{n-1}$ in $\R^n$. We are 
looking for natural general conditions on $X_k$ which guarantee 
that the distribution functions $F_\theta(x) = \P\{S_\theta \leq x\}$
are well approximated for most of $\theta \in \S^{n-1}$ by the standard 
normal distribution function
$$
\Phi(x) = \frac{1}{\sqrt{2\pi}} \int_{-\infty}^x e^{-y^2/2}\,dy, 
\qquad x \in \R.
$$
Of special interest is the question of possible rates
in the Kolmogorov distance
$$
\rho(F_\theta,\Phi) = \sup_x |F_\theta(x) - \Phi(x)|.
$$

In this problem, going back to the seminal work of Sudakov \cite{Su},
the well studied classical case of independent components may serve 
as a basic example for comparison with various models or dependencies. 
Let us recall that, if $X_k$ are independent and have finite 4-th moments 
(with mean zero and variance one), there is an upper bound on average
\be
c\,\E_\theta\, \rho(F_\theta,\Phi) \leq \frac{1}{n}\,\bar \beta_4, \quad
\bar \beta_4 = \frac{1}{n} \sum_{k=1}^n \E X_k^4,
\en
where $c>0$ is an absolute constant, and where we use $\E_\theta$ 
to denote the integral over the uniform probability measure $\s_{n-1}$ 
on the unit sphere. Moreover, for any $r>0$,
\be
\s_{n-1}\Big\{c\,\rho(F_\theta,\Phi) \geq \frac{1}{n}\,
\bar \beta_4 \, r\Big\} \, \leq \, 2\,e^{-\sqrt{r}}.
\en
This non-trivial phenomenon was observed by Klartag and Sodin \cite{K-S}. 
It shows that when $\bar \beta_4$ is bounded like in the i.i.d. situation, 
the distances $\rho(F_\theta,\Phi)$ turn out to be typically of order 
at most $1/n$. This is in contrast to the case of equal coefficients leading to 
the unimprovable standard $\frac{1}{\sqrt{n}}$-rate (in general, including
independent Bernoulli summands $X_k$). Moreover, in the i.i.d. situation with 
finite moment $\beta_5 = \E\, |X_1|^5$ and symmetric underlying distributions, 
the typical rate of normal approximation for $F_\theta$ may further be improved 
to $\beta_5\, n^{-3/2}$ up to a constant (which is best possible as long as 
$\E X_1^4 \neq 3$, cf. \cite{B5}).

As for more general models with not necessarily independent components $X_k$, 
the study of this high-dimensional phenomenon has a long history, and we refer 
an interested reader to the book \cite{B-G-V-V} and a recent paper \cite{B-C-G5} 
for an account of various results in this direction. Let us only mention 
\cite{A-B-P}, \cite{B2}, \cite{B3}, \cite{So}, \cite{K1}, \cite{K2}, \cite{E-K},
where one can find quantitative variants of Sudakov's theorem on the 
concentration of $F_\theta$ about the typical (average) distribution 
$F = \E_\theta F$ and/or about the normal law $\Phi$ for different metrics and 
under certain assumptions (of convexity-type, for example). Some papers provide 
Berry-Esseen-type estimates on the closeness of $F_\theta$ to $\Phi$ 
explicitly in terms of $\theta$ assuming that the distribution of the random 
vector $X$ is ``sufficiently" symmetric, cf. \cite{M-M}, \cite{M}, \cite{G-S}, 
\cite{K3}, \cite{G}.

Whether or not $F$ itself is close to the standard normal law represents 
a thin-shell problem on the concentration of the values of the square 
of the Euclidean norm $|X|$ about its mean $\E\,|X|^2 = n$ (or, in essence, on the 
concentration of $|X|$ about $\sqrt{n}$). The rate of concentration may be 
controlled in terms of the functional $\sigma_4 = \frac{1}{n}\,\Var(|X|^2)$ 
which is often of order 1 (including the i.i.d. situation). Once it is the case,
one can obtain a standard rate of concentration of $F_\theta$ around $\Phi$
on average under mild moment assumptions. For example, it is known that, 
if $\E\,|X|^2 = n$ (without the isotropy hypothesis), then 
$$
\E_\theta\, \rho(F_\theta,\Phi) \, \leq \,
c\,\big(M_3^3 + \sigma_4^{3/2}\big) \frac{1}{\sqrt{n}},
$$
up to an absolute constant $c>0$,
where $M_3^3 = \sup_{\theta} \E\,|S_\theta|^3$ (cf. \cite{B-C-G4}).

In order to reach better rates, one has to involve
stronger assumptions or functionals such as $\Lambda = \Lambda(X)$ defined
as an optimal constant in the inequality
\be
\Var\bigg(\sum_{i,j=1}^n a_{ij} X_i X_j\bigg) \leq 
\Lambda \sum_{i,j=1}^n a_{ij}^2 \qquad (a_{ij} \in \R),
\en
which may be referred to as a second order correlation condition.
In terms of $\Lambda$, the bound (1.1) has been extended in \cite{B-C-G5} modulo 
a logarithmic factor: If additionally to the isotropy assumption the distribution 
of $X$ is symmetric around the origin, it was shown that
\be
c\,\E_\theta\, \rho(F_\theta,\Phi) \leq \frac{\log n}{n}\, \Lambda.
\en

The optimal value $\Lambda = \Lambda(X)$ in (1.3) is finite as long 
as $|X|$ has a finite $4$-th moment. It represents the maximal 
eigenvalue of the covariance matrix associated to the 
$n^2$-dimensional random vector $(X_i X_j - \E X_i X_j)_{i,j=1}^n$.
This parameter may be effectively estimated in many examples
and is related to other standard characteristics. For example, 
$\Lambda(X) \leq 2\,\max_k\,\E X_k^4$, if $X_k$ are independent. If 
$X$ is isotropic, and its distribution admits a Poincar\'e-type inequality
\be
\lambda_1\, \Var(u(X)) \, \leq \, \E\, |\nabla u(X)|^2
\en
with a positive (optimal) constant $\lambda_1 = \lambda_1(X)$
for all smooth functions $u$ on $\R^n$, then we have 
$\Lambda(X) \leq 4/\lambda_1(X)$.

The aim of these notes is to sharpen (1.4) via a large deviation 
bound in analogy with (1.2). This turns out to be possible
as long as all linear forms $S_\theta$ have finite exponential
moments. To avoid technical discussions, we restrict ourselves 
to the case where $\lambda_1 > 0$, which at the same time allows
to drop the symmetry assumption.

\vskip5mm
{\bf Theorem 1.1.} {\sl Let $X$ be an isotropic random vector 
in $\R^n$ with mean zero and a positive Poincar\'e constant 
$\lambda_1$. Then with some absolute constant $c>0$
\be
c\,\E_\theta\, \rho(F_\theta,\Phi) \, \leq \, 
\frac{\log n}{n}\, \lambda_1^{-1}.
\en
Moreover, for all $r > 0$,
\be
\s_{n-1}\Big\{c\,\rho(F_\theta,\Phi) \geq 
\frac{\log n}{n}\,\lambda_1^{-1} \, r\Big\} \, \leq \, 
2\,e^{-\sqrt{r}}.
\en
}

\vskip2mm
Being restricted to isotropic log-concave distributions, 
an interesting feature of the bound (1.4) is its connection
with certain open problems in Asymptotic Convex Geometry
such as the K-L-S and thin-shell conjectures.
Namely, modulo $n$-dependent logarithmic factors, the following 
three assertions are equivalent up to positive constants
$c$ and $\beta$ (perhaps different in different places) for the
entire class of isotropic random vectors $X$ in $\R^n$
having symmetric log-concave distributions (cf. \cite{B-C-G5}):

\vskip5mm
(i) $\sup_X \lambda_1^{-1}(X) \leq c\,(\log n)^\beta$;
\qquad \ \ 
(ii) $\sup_X \Var(|X|) \leq c\,(\log n)^\beta$;

\vskip2mm
(iii) $\sup_X \E_\theta\, 
\rho(F_\theta,\Phi) \leq \frac{c}{n}\,(\log n)^\beta$.

\vskip5mm 
\noindent
In this connection, let us also mention a recent paper by Jiang,
Lee and Vempala \cite{J-L-V}, which provides a reformulation
of (i)-(ii) as a central limit theorem for random variables
of the form $\left<X,Y\right>$, where $Y$ is an independent copy of $X$.

Note that the implication ${\rm (i)} \Rightarrow {\rm (ii)}$ is immediate
when applying (1.5) to $u(x) = |x|$, while the reverse statement 
is a deep theorem due to Eldan \cite{E}. By (1.4), we also have 
${\rm (i)} \Rightarrow {\rm (iii)}$. As for the implication
${\rm (iii)} \Rightarrow {\rm (ii)}$, it holds true in view of 
a general relation
$$
c\,\Var(|X|) \, \leq \,
n\,(\log n)^4\ \E_\theta\, \rho(F_\theta,\Phi) + 1
$$
(which only requires that all $S_\theta$ have a finite and bounded
exponential moment).

The symmetry assumption is irrelevant both in (i) and (ii). However, this is 
not so obvious concerning (iii). Indeed, one may try to use a symmetrization 
argument by applying (1.4) to the random vector $X' = (X - Y)/\sqrt{2}$. 
But then we need a quantitative form of a particular variant of Cramer's 
theorem: If $\eta$ is an independent copy of a random variable $\xi$ with 
mean zero and variance one, and if $\xi' = (\xi - \eta)/\sqrt{2}$ is almost 
standard normal, then so is $\xi$. The best result in this direction 
is the following theorem due to Sapogov \cite{Sa}: Given that 
$\rho(F_{\xi'},\Phi) \leq \ep \leq 1/e$, we have
$$
\rho(F_\xi,\Phi) \leq C\,\big(\log(1/\ep)\big)^{-1/2}
$$
up to some absolute constant $C$, where $F_\xi$ and $F_{\xi'}$ 
denote the distribution functions of $\xi$ and $\xi'$.
Moreover, the dependence in $\ep$ on the right-hand side
cannot be improved, as was shown in \cite{C} (cf. also
\cite{B-C-G1} for a related model).
Thus, the resulting bound on $\E_\theta\, \rho(F_\theta,\Phi)$ 
which can be derived this way on the basis of $X'$ cannot
yield even a standard rate.

Here, we choose a different route. As we will see, it is possible
to remove the symmetry hypothesis, by adding to the right-hand 
side of (1.4) an additional term responsible for higher order 
correlations between $X_k$. More precisely, as a preliminary bound
which is based on the $\Lambda$-functional only, it will be shown that
\be
c\,\E_\theta\, \rho(F_\theta,\Phi) \, \leq \,
\frac{\log n}{n}\, \Lambda + \Big(\frac{\log n}{n}\Big)^{1/4}
\bigg(\E\, \frac{\left<X,Y\right>}{\sqrt{|X|^2 + |Y|^2}}\bigg)^{1/2}.
\en
The last expectation is vanishing for symmetric distributions, or, 
for example, if $|X| = \sqrt{n}$ a.s. As another scenario, 
the second term in (1.8) is of a smaller order in comparison with 
$\frac{\log n}{n}\, \lambda_1^{-1}$ when (1.5) holds. Nevertheless, 
in contrast to the bound (1.4), the derivation of (1.8) turns out 
to be tedious, since it involves a careful analysis of 
projections of the characteristic functions $f_\theta(t)$ of 
$S_\theta$ as functions of $\theta$ onto the subspace of all 
linear functions in the Hilbert space $L^2(\R^n,\s_{n-1})$.

The paper is organized as follows. We start with the study of
densities of linear functionals on the sphere $\S^{n-1}$ viewed 
as random variables with respect to the normalized Lebesgue measure 
$\s_{n-1}$. Here, the aim will be to refine the asymptotic normality 
of these distributions in analogy with Erdgeworth expansions in the 
central limit theorem (which we consider up to order 2, Sections 2-3). 
Then we turn to the problem of deviations of general smooth functions 
on $\S^{n-1}$ in terms of their Hessians, recalling and extending 
several results in this direction (Section~4). These results are 
applied in Sections 5 to characteristic functions $f_\theta(t)$, 
with a separate treatment of their linear parts in $L^2(\s_{n-1})$ 
in the next Section 6. In Section 7, we adapt basic Fourier 
analytic tools in the form of Berry-Esseen-type bounds to the 
scheme of weighted sums. Deviations of involved integrals as 
functions on the sphere are discussed separately in Section 8.
Section 9 collects several general facts about Poincar\'e-type
inequalities that will be needed for the proof of Theorem 1.1, while
final steps of the proof are deferred to the remaining Sections 10-12.

As usual, the Euclidean space $\R^n$ is endowed with the canonical norm 
$|\,\cdot\,|$ and the inner product $\left<\cdot,\cdot\right>$. 
We denote by $c$ a positive absolute constant which may vary
from place to place (if not stated explicitly that $c$ depends
on some parameter).


\vskip7mm
\section{{\bf Distribution of Linear Functionals on the Sphere}}
\setcounter{equation}{0}

\vskip2mm
\noindent
By the rotational invariance of $\s_{n-1}$, all linear functionals 
$u(\theta) = \left<\theta,v\right>$ with $|v|=1$ have equal 
distributions. Hence, it is sufficient to focus just on the first coordinate 
$\theta_1$ of the vector $\theta \in \S^{n-1}$ viewed as a random variable
on the probability space $(\S^{n-1},\s_{n-1})$. 
It is well-known that this random variable has density
$$
c_n\, \big(1 - x^2\big)_+^{\frac{n-3}{2}}, \quad x \in \R, \qquad
c_n =  \frac{\Gamma(\frac{n}{2})}{\sqrt{\pi}\,\Gamma(\frac{n-1}{2})},
$$
with respect to the Lebesgue measure on the real line, where 
$c_n$ is a normalizing constant. 

We denote by $\varphi_n$ the density of the normalized first 
coordinate $\sqrt{n}\, \theta_1$, i.e.,
$$
\varphi_n(x) =  c_n'\, \bigg(1 - \frac{x^2}{n}\bigg)_+^{\frac{n-3}{2}}, 
\quad c_n' = \frac{c_n}{\sqrt{n}}.
$$
Clearly, 
$$
\varphi_n(x) \rightarrow \varphi(x) = \frac{1}{\sqrt{2\pi}}\,e^{-x^2/2},
\qquad
c_n'  \rightarrow \frac{1}{\sqrt{2\pi}} = 0.399...
$$
as $n \rightarrow \infty$, and one can also show that $
c_n'  < \frac{1}{\sqrt{2\pi}}$ for all $n$.

Deviations for $\varphi_n(x)$ from $\varphi(x)$ have been
considered in \cite{B-C-G4}. In particular, if $n \geq 3$, then
for all $x \in \R$,
\be
|\varphi_n(x) - \varphi(x)| \, \leq \, \frac{c}{n}\,e^{-x^2/4}.
\en
We need to sharpen this bound by obtaining an approximation for 
$\varphi_n(x)$ with an error of order $1/n^2$ by means of 
a suitable modification of the standard normal density. Denote by 
$H_4(x) = x^4 - 6x^2 + 3$ the 4-th Chebyshev-Hermite polynomial.

\vskip5mm
{\bf Proposition 2.1.} {\sl For all $x \in \R$ and $n \geq 3$, 
\be
\Big|\varphi_n(x) - \varphi(x) \Big(1 - \frac{H_4(x)}{4n}\Big)\Big|
 \, \leq \, \frac{c}{n^2}\,e^{-x^2/4}.
\en
}

\vskip2mm
{\bf Proof.} In the interval $|x| \leq \frac{1}{2}\sqrt{n}$, 
consider the function $p_n(x) = (1 - \frac{x^2}{n})_+^{\frac{n-3}{2}}$. 
Using the Taylor expansion for the logarithmic function near zero, 
one may write
\bee
-\log p_n(x)
 & = &
-\frac{n-3}{2}\,\log\bigg(1 - \frac{x^2}{n}\bigg) \\
 & = &
\frac{n-3}{2}\,\left(\frac{x^2}{n} + \frac{x^4}{2n^2} + 
\bigg(\frac{x^2}{n}\bigg)^3 \sum_{k=3}^\infty \frac{1}{k}\, 
\bigg(\frac{x^2}{n}\bigg)^{k-3}\right)
 \ = \
\frac{x^2}{2} + \delta.
\ene
The remainder term has the form 
$$
\delta =  - \frac{3x^2}{2n} + \frac{x^4}{4n} - \frac{1}{n^2}\,
\Big(\frac{3}{4}\,x^4 - \frac{n-3}{3n}\, x^6 \ep\Big)
$$
with some $0 \leq \ep \leq 1$. By the assumption that
$x^2 \leq \frac{1}{4}\,n$, it satisfies
$$
\delta 
 \, \geq \,
- \frac{3x^2}{2n} + \frac{x^4}{4n} - \frac{3x^4}{4n^2}
 \, \geq \,
- \frac{27x^2}{16n} + \frac{x^4}{4n} \ > \ -\frac{27}{64}.
$$
Hence
$$
|e^{-\delta} - 1 + \delta| \leq \frac{\delta^2}{2}\,e^{27/64} \leq 
\delta^2.
$$
Moreover, using once more $x^2 \leq \frac{1}{4}\,n$, we get 
$$
|\delta| 
 \, \leq \,
\frac{3x^2}{2n} + \frac{x^4}{4n} + \frac{1}{n^2}\,
\bigg(\frac{3}{4}\,x^4 + \frac{1}{3}\, x^6\bigg)
 \, \leq \,
\frac{x^2}{n}\, \bigg(\frac{27}{16} + \frac{1}{3}\, x^2\bigg), \\
$$
which implies
$$
\delta^2 \leq \frac{x^4}{n^2}\, \bigg(6 + \frac{2}{9}\, x^4\bigg).
$$
Hence, with some $|\ep_1| \leq 1$,
$$
e^{x^2/2}\,p_n(x) 
 \, = \,
e^{-\delta} \, = \, 1 - \delta + \ep_1 \delta^2
 \, = \,
1 +  \frac{3x^2}{2n} - \frac{x^4}{4n} + \frac{A}{n^2},
$$
where 
\bee
|A| 
 & \leq &
\bigg|\,\frac{3}{4}\,x^4 - \frac{n-3}{3n}\, x^6 \ep\bigg| + 
x^4\, \bigg(6 + \frac{2}{9}\, x^4\bigg) \\
 & \leq &
\frac{3}{4}\,x^4 + \frac{1}{3}\, x^6 + 
x^4\, \bigg(6 + \frac{2}{9}\, x^4\bigg) 
 \ \leq \
8x^4 + x^8.
\ene
As a result,
\be
p_n(x) = e^{-x^2/2}\, 
\bigg[1 + \frac{6x^2 - x^4}{4n} + \frac{\ep}{n^2}\, 
\big(8 x^4 + x^8\big)\bigg], \qquad |\ep| \leq 1.
\en

To derive a similar expansion for $\varphi_n(x)$, denote by 
$Z$ a standard normal random variable. From (2.3) we obtain that
\bee
\frac{1}{\sqrt{2\pi}} \int_{-\infty}^\infty p_n(x)\,dx 
 & = & 
1 + \frac{1}{4n}\,\big(6\, \E Z^2 - \E Z^4\big) + 
O\Big(\frac{1}{n^2}\Big) \\
 & = & 
1 + \frac{3}{4n} + O\Big(\frac{1}{n^2}\Big).
\ene
Here we used the property that $p_n(x)$ has a sufficiently fast decay
for $|x| \geq \frac{1}{2}\sqrt{n}$, as indicated in (2.1). Since 
$\varphi_n(x) = c_n'\,p_n(x)$ is a density, we conclude that
$$
1 = c_n'
\sqrt{2\pi}\,\bigg(1 + \frac{3}{4n} + O\Big(\frac{1}{n^2}\Big)\bigg),
\qquad
c_n'\sqrt{2\pi} = 1 - \frac{3}{4n} + O\Big(\frac{1}{n^2}\Big).
$$
Hence
\bee
\sqrt{2\pi}\ e^{x^2/2}\, \varphi_n(x)
 & = & 
\bigg(1 - \frac{3}{4n} + O\Big(\frac{1}{n^2}\Big)\bigg) \, 
\bigg[1 + \frac{6x^2 - x^4}{4n} + 
\frac{\ep}{n^2}\, \Big(8x^4 + x^8\Big)\bigg] \\
 & = &
1 + \frac{6x^2 - x^4 }{4n} - \frac{3}{4n} + 
O\bigg(\frac{1 + x^8}{n^2}\bigg).
\ene
Thus, in the interval $|x| \leq \frac{1}{2}\sqrt{n}$,
$$
\varphi_n(x) = \varphi(x)\, \bigg[1 - \frac{H_4(x)}{4n} + Q_n(x)\,
\frac{1 + x^8}{n^2}\bigg]
$$
with a quantity $Q_n(x)$ bounded by a universal constant
in absolute value. In view of (2.1), the bound (2.2) follows immediately.
\qed


\vskip7mm
\section{{\bf Characteristic Function of Linear Functionals}}
\setcounter{equation}{0}

\vskip2mm
\noindent
In the sequel, we denote by $J_n = J_n(t)$ the characteristic function 
of the first coordinate $\theta_1$ of a random vector 
$\theta = (\theta_1,\dots,\theta_n)$ which is 
uniformly distributed on the unit sphere $\S^{n-1}$.
In a more explicit form, for any $t \in \R$,
$$
J_n(t) 
 \, = \,
c_n \int_{-\infty}^{\infty} e^{itx}\,(1 - x^2)_+^{\frac{n-3}{2}}\,dx 
 \, = \,
c_n' \int_{-\infty}^{\infty} e^{itx/\sqrt{n}}\,
\Big(1 - \frac{x^2}{n}\Big)_+^{\frac{n-3}{2}}\,dx.
$$
This is just a multiple of the Bessel function 
of the first kind with index $\nu = \frac{n}{2} - 1$
(\cite{Bat}, p.\,81).

Thus, the characteristic function of the 
normalized first coordinate $\theta_1 \sqrt{n}$ is given by
$$
\hat \varphi_n(t) = J_n\big(t\sqrt{n}\big) = 
\int_{-\infty}^\infty e^{itx}\varphi_n(x)\,dx,
$$
which is the Fourier transform of the probability density $\varphi_n$. 
Proposition 2.1 can be used to compare $\hat \varphi_n(t)$ with the Fourier 
transform of the ``corrected Gaussian measure'', as well as to compare 
the derivatives of these transforms. 

\vskip5mm
{\bf Proposition 3.1.} {\sl For all $t \in \R$,
$$
\bigg|J_n\big(t\sqrt{n}\big) - \Big(1 - \frac{t^4}{4n}\Big)\, 
e^{-t^2/2}\bigg| \, \leq \, \frac{c}{n^2}.
$$
Moreover, for any $k = 1,2,\dots$,
$$
\bigg|\frac{d^k}{dt^k}\, J_n\big(t\sqrt{n}\big) - \frac{d^k}{dt^k}\, 
\bigg(\Big(1 - \frac{t^4}{4n}\Big)\, e^{-t^2/2}\bigg)\bigg|
 \, \leq \, \frac{(ck)^{k/2}}{n^2}.
$$
}

\vskip4mm
Taking $k=1$, we have
$$
\bigg|\big(J_n\big(t\sqrt{n}\big)\big)' - 
\Big(\frac{t^5}{4n} - \frac{t^3}{n} - t\Big)\, e^{-t^2/2}\bigg|
 \, \leq \, \frac{c}{n^2}.
$$
One may also add a $t$-depending factor on the right-hand side. 
For $t$ of order 1, this can be done just by virtue of Taylor's 
formula. Indeed, the functions
$$
f_n(t) = J_n\big(t\sqrt{n}) = \E_\theta\,e^{it\theta_1 \sqrt{n}}, \qquad
g_n(t) = \Big(1 - \frac{t^4}{4n}\Big)\, e^{-t^2/2}
$$
have equal first three derivatives at zero. Since, by Proposition 3.1,
$|f_n^{(4)}(t) - g_n^{(4)}(t)| \leq \frac{c}{n^2}$, Taylor's 
formula refines this proposition for the interval $|t| \leq 1$.

\vskip5mm
{\bf Corollary 3.2.} {\sl For all $t \in \R$,
$$
\Big|J_n\big(t\sqrt{n}\big) - \Big(1 - \frac{t^4}{4n}\Big)\, 
e^{-t^2/2}\Big| \, \leq \, \frac{c}{n^2}\,\min\{1,t^4\},
$$
$$
\Big|\big(J_n\big(t\sqrt{n}\big)\big)' +
t\,\Big(1 + \frac{4t^2 - t^4}{4n}\Big)\, e^{-t^2/2}\Big|
 \, \leq \, \frac{c}{n^2}\,\min\{1,|t|^3\}.
$$
}

\vskip5mm
These approximations may be complemented by a Gaussian decay bound 
\be
\big|J_n\big(t\sqrt{n}\big)\big| \, \leq \,
5\,e^{-t^2/2} + 4\,e^{-n/12}, \quad t \in \R
\en
(cf. \cite{B-C-G4}, Proposition 3.3).

\vskip5mm
{\bf Proof of Proposition 3.1.} In general, given two integrable 
functions on the real line, say, $p$ and $q$, their Fourier transforms
$$
\hat p(t) = \int_{-\infty}^\infty e^{itx} p(x)\,dx, \qquad
\hat q(t) = \int_{-\infty}^\infty e^{itx} q(x)\,dx
$$
satisfy, for all $t \in \R$, 
$$
|\hat p(t) - \hat q(t)| \leq \int_{-\infty}^\infty|p(x) - q(x)|\,dx.
$$
Moreover, one may differentiate these transforms $k$ times to get
$$
\frac{d^k}{dt^k}\, \hat p(t) \, = \,
\int_{-\infty}^\infty (ix)^k\,e^{itx}\,p(x)\,dx, \qquad
\frac{d^k}{dt^k}\, \hat q(t) \, = \,
\int_{-\infty}^\infty (ix)^k\,e^{itx}\,q(x)\,dx,
$$
as long as the integrands are integrable, which also yields the relation
$$
\Big|\frac{d^k}{dt^k}\, \hat p(t) - \frac{d^k}{dt^k}\, \hat q(t)\Big| \, \leq \,
\int_{-\infty}^\infty |x|^k\,|p(x) - q(x)|\,dx.
$$
This applies in particular to the functions $p(x) = \varphi_n(x)$ 
and $q(x) = \varphi(x)\,(1 - \frac{1}{4n}\,H_4(x))$ whose
Fourier transform is described as
$$
\hat q(t) = e^{-t^2/2} \Big(1 - \frac{t^4}{4n}\Big).
$$
Since (by Stirling's formula)
$$
\int_{-\infty}^\infty |x|^k\,e^{-x^2/4}\,dx \, = \,
2^{k+1}\,\Gamma\Big(\frac{k+1}{2}\Big) \, \leq \, (ck)^{k/2},
$$
it remains to apply (2.2).
\qed


\vskip7mm
\section{{\bf Deviations of Smooth Functions on the Sphere}}
\setcounter{equation}{0}

\vskip2mm
\noindent
Smooth functions $u$ on the unit $n$-sphere with $\s_{n-1}$-mean zero
are known to have fluctuations of order at most $1/\sqrt{n}$
(which is the case for all linear functions).
This may be seen from the Poincar\'e inequality
\be
\int |u|^2\,d\s_{n-1} \leq \frac{1}{n-1}\int |\nabla u|^2\,d\s_{n-1}.
\en
Moreover, when $u$ is Lipschitz, that is, $|\nabla u(\theta)| \leq 1$
for all $\theta \in \S^{n-1}$, there is a subgaussian exponential
bound on the Laplace transform (cf. \cite{L2})
\be
\int \exp\Big\{\sqrt{n-1}\, ru\Big\}\,d\s_{n-1} \, \leq \, 
e^{r^2/2}, \qquad r \in \R.
\en

This spherical concentration phenomenon may be strengthened with 
respect to the dimension $n$ for a wide subclass of smooth functions. 
We denote by $\nabla^2 u(x)$ the Hessian, that is, the $n \times n$ 
matrix of second order partial derivative $\partial_{ij} u(x)$, and by
$\I_n$ the identity $n \times n$ matrix. Given a symmetric matrix 
$A = (a_{ij})_{i,j=1}^n$ with real or complex entries, 
the associated Hilbert-Schmidt and operator norms are defined by 
$$
\|A\|_{\rm HS} = \bigg(\sum_{i,j=1}^n |a_{ij}|^2\bigg)^{1/2}, \qquad
\|A\| = \max_{|\theta| = 1} |\left<A\theta,\theta\right>|.
$$

The next proposition summarizes several results from \cite{B-C-G5}
employing a second order concentration on the sphere, a property 
developed in \cite{B-C-G2}.

\vskip5mm
{\bf Proposition 4.1.} {\sl Suppose that a real-valued function $u$ is 
defined and $C^2$-smooth in some neighborhood of $\S^{n-1}$.
If $u$ is orthogonal to all affine functions in $L^2(\s_{n-1})$, then
\be
\int |u|^2\,d\s_{n-1} \leq \frac{5}{(n-1)^2}
\int \|\nabla^2 u - a\,\I_n\|_{\rm HS}^2\,d\s_{n-1}
\en
for any $a \in \R$. Moreover, if $\|\nabla^2 u - a\, \I_n\| \leq 1$
on $\S^{n-1}$ and the second integral in $(4.3)$ is bounded by $b$, then
\be
\int \exp\Big\{\frac{n-1}{2(1+4b)}\, |u|\Big\}\, d\s_{n-1}
 \, \leq \, 2.
\en
}

\vskip5mm
By Markov's inequality, (4.4) yields a corresponding large deviation 
bound, which may be stated informally as a subexponential stochastic dominance
$|u| \leq c_b\,(\frac{1}{\sqrt{n}} Z)^2$ with $Z \sim N(0,1)$. 
Thus, the deviations of $u$ are of order at most $1/n$.

We will need the following generalization of Proposition 4.1 which is more 
flexible in applications. Given a function $u$ in the (complex) Hilbert space 
$L^2 = L^2(\R^n,\s_{n-1})$, we consider its orthogonal projection
$$
l = {\rm Proj}_H u
$$
onto the linear space $H$ in $L^2$ generated by the constant and linear functions
on $\R^n$. Let us call $l$ an affine part of $u$.

\vskip5mm
{\bf Proposition 4.2.} {\sl Suppose that a complex-valued function 
$u$ is $C^2$-smooth in some neighborhood of $\S^{n-1}$ and has 
$\s_{n-1}$-mean zero. For any $a \in \C$, 
\be
\int |u|^2\,d\s_{n-1} \, \leq \, \frac{5}{(n-1)^2}
\int \|\nabla^2 u - a\,\I_n\|_{\rm HS}^2\,d\s_{n-1} + \|l\|_{L^2}^2,
\en
where $l$ is the affine part of $u$. Moreover, if 
$\|\nabla^2 u - a\, \I_n\| \leq 1$ on $\S^{n-1}$, then
\be
\|u\|_{\psi_1} \, \leq \, \frac{4}{n-1} + \frac{16}{n-1}
\int \|\nabla^2 u - a\,\I_n\|_{\rm HS}^2\,d\s_{n-1} + 6\,\|l\|_{L^2}.
\en
}

\vskip5mm
Here we used a standard notation
$$
\|u\|_{\psi_1} = \inf\Big\{\lambda > 0: \E_\theta\,e^{|u|/\lambda} \leq 2\Big\}
$$
for the Orlicz norm on the probability space $(\S^{n-1},\s_{n-1})$
generated by the Young function $\psi_1(r) = e^{|r|} - 1$ ($r \in \R$).

\vskip5mm
{\bf Proof of Proposition 4.2.}
The Poincar\'e-type inequalities (4.1) and (4.3) continue to hold
in the class of all complex-valued functions $u$ with $\s_{n-1}$-mean zero,
while (4.2) and (4.4) require slight modifications.
Indeed, (4.4) may be applied separately to the real part
$u_0 = {\rm Re}(u)$ and to the imaginary part $u_1 = {\rm Re}(u)$ of $u$, 
which results in
\be
\int \exp\Big\{\frac{n-1}{2(1+2b_k)}\, |u_k|\Big\}\, d\s_{n-1} \leq 2, \qquad
b_k = \int \|\nabla^2 u_k - a_k\,\I_n\|_{\rm HS}^2\,d\s_{n-1},
\en
for $k=0$ and $k=1$, assuming that the following conditions are fulfilled:

\vskip4mm
$a)$\, $u_0$ and $u_1$ (that is, $u$) are $C^2$-smooth and orthogonal to all 
affine functions in $L^2(\s_{n-1})$; 

$b)$\, $\|\nabla^2 u_k - a_k\, \I_n\| \leq 1$ on $\S^{n-1}$ with
$a_0 = {\rm Re}(a)$ and $a_1 = {\rm Re}(a)$.

\vskip4mm
\noindent
The latter requirement is met as long as
\be
\|\nabla^2 u - a\, \I_n\| \, \equiv \,
\max_{|\theta| = 1} |\left<(\nabla^2 u - a\, \I_n)\theta,\theta\right>| \,\leq \, 1
\en
pointwise on $\S^{n-1}$. As for the exponential bounds in (4.7), they may 
equivalently be written in terms of the Orlicz $\psi_1$-norm as
$$
\|u_k\|_{\psi_1} \, \leq \, \frac{2}{n-1} + \frac{8b_k}{n-1}, \quad k = 0,1.
$$
Applying the triangle inequality 
$\|u\|_{\psi_1} \leq \|u_0\|_{\psi_1} + \|u_1\|_{\psi_1}$ in the Orlicz space
and noting that $b_0 + b_1$ is just the integral on the right-hand side
in (4.5)-(4.6), we conclude that
\be
\|u\|_{\psi_1} \, \leq \, \frac{4}{n-1} + \frac{16}{n-1}
\int \|\nabla^2 u - aI_n\|_{\rm HS}^2\,d\s_{n-1}.
\en
This is a ``complex" variant of the inequality (4.4), which holds for all
$a \in \mathbb C$ under the assumption that $u$ is $C^2$-smooth in some 
neighborhood of $\S^{n-1}$, is orthogonal to all affine functions in 
$L^2(\s_{n-1})$, and satisfies (4.8). 

One may now start with an arbitrary $C^2$-smooth function $u$ with mean zero, 
but apply these hypotheses and the conclusions to the 
projection $Tu$ of $u$ onto the orthogonal complement of the space $H$ 
of all linear functions in $L^2(\s_{n-1})$. This space has dimension $n$, 
and one may choose for the orthonormal basis in $H$ the canonical functions 
$$
l_k(\theta) = \sqrt{n}\,\theta_k, \ \ 
k = 1,\dots,n, \ \ \theta = (\theta_1,\dots,\theta_n) \in \S^{n-1}.
$$
Therefore, the ``linear" part $l = Tu - u$ of $u$ is described 
as the orthogonal projection in $L^2(\s_{n-1})$ onto $H$, namely
\bee
l(\theta)
 & = &
\sum_{k=1}^n \left<u,l_k\right>_{L^2} l_k(\theta)
 \ = \
\sum_{k=1}^n \Big(\int u(x) l_k(x)\,d\s_{n-1}(x)\Big)\, l_k(\theta) \\
 & = & 
n \int \Big(u(x)\,\sum_{k=1}^n x_k \theta_k\Big)\,d\s_{n-1}(x).
\ene
In other words, 
$$
l(\theta) = \left<v,\theta\right> \quad {\rm with} \ 
v = n \int xu(x)\,d\s_{n-1}(x),
$$
which implies, in particular, that
\be
\|l\|_{L^2}^2 = \frac{1}{n}\,|v|^2 = n
\int\!\!\! \int
\left<x,y\right> u(x)\bar u(y)\,d\s_{n-1}(x)d\s_{n-1}(y).
\en

The functions $Tu$ and $u$ have identical Euclidean second derivatives. 
Hence, (4.5) follows from (4.3) when the latter is applied to $Tu$,
since $Tu$ and $l$ are orthogonal in $L^2$. Applying (4.9) with $Tu$ 
in place of $u$, we similarly have
\be
\|Tu\|_{\psi_1} \, \leq \, \frac{4}{n-1} + \frac{16}{n-1}
\int \|\nabla^2 u - aI_n\|_{\rm HS}^2\,d\s_{n-1},
\en
provided that
$\|\nabla^2 Tu - a\, \I_n\| = \|\nabla^2 u - a\, \I_n\| \leq 1$ on $\S^{n-1}$
as in (4.8).

To derive (4.6), it remains to use the fact that the linear functions on the sphere 
behave like Gaussian random variables. This can be seen from (4.2), which may
be applied with $r=1$ to the real and imaginary parts
of $l/\|l\|_{{\rm Lip}}$. Then it gives
$$
\int \exp\Big\{\sqrt{n-1}\,|l|/4\,\|l\|_{{\rm Lip}}\Big\}\,d\s_{n-1} \leq 2,
$$
so that
$$
\|l\|_{\psi_1} \leq \frac{4}{\sqrt{n-1}}\,\|l\|_{{\rm Lip}} =
\frac{4\sqrt{n}}{\sqrt{n-1}}\,\|l\|_{L^2} \leq 6\,\|l\|_{L^2}.
$$
The latter should be combined with (4.11), and we arrive at (4.6) due to the
triangle inequality $\|u\|_{\psi_1} \leq \|Tu\|_{\psi_1} + \|l\|_{\psi_1}$.
\qed


\vskip7mm
\section{{\bf Concentration of Characteristic Functions}}
\setcounter{equation}{0}

\vskip2mm
\noindent
Given a random vector $X = (X_1,\dots,X_n)$ in $\R^n$, we consider
the smooth functions
\be
u_t(\theta) = f_\theta(t) = \E\,e^{it\left<X,\theta\right>}, 
\quad \theta \in \R^n,
\en
where $t \in \R$ serves as a parameter. For any fixed $\theta \in \R^n$, 
$t \rightarrow f_\theta(t)$ represents the characteristic function 
of the weighted sum $S_\theta = \left<X,\theta\right>$ with 
distribution function $F_\theta$, while the $\s_{n-1}$-mean
$$
f(t) = \E_\theta f_\theta(t) = \E_\theta\, \E\,e^{it\left<X,\theta\right>}
$$
is the characteristic function of the average distribution function
$F(x) = \E_\theta\, \P\{S_\theta \leq x\}$, $x \in \R$. 
Recall that we use $\E_\theta$ to denote integrals
with respect to the uniform measure $\s_{n-1}$.

In order to control deviations of $u_t$ from $f(t)$ on $\S^{n-1}$ at 
the standard rate, the spherical concentration inequalities (4.1)-(4.2) 
are sufficient. Indeed, the function $u_t$ has a gradient described 
in the vector form as
$$
\left<\nabla u_t(\theta),w\right> \, = \, 
it\, \E \left<X,w\right> e^{it \left<X,\theta\right>}, 
\qquad w \in \C^n.
$$
Hence, under the isotropy assumption, writing $w = w_0 + iw_1$
($w_0,w_1 \in \R^n$), we have
\bee
|\left<\nabla u_t(\theta),w\right>|^2 
 & \leq &
\E\, |\left<X,w\right>|^2 \\
 & = &
\E \left<X,w_0\right>^2 + \E \left<X,w_1\right>^2 \, = \,
|w_0|^2 + |w_1|^2 \, = \, |w|^2,
\ene
that is, $|\left<\nabla u_t(\theta),w\right>| \leq |t|\,|w|$ for all 
$w \in \C^n$. This gives a uniform bound $|\nabla u_t(\theta)| \leq |t|$,
so that, by the spherical Poincar\'e inequality (4.1),
\be
\E_\theta\,|f_\theta(t) - f(t)|^2 \, \leq \, \frac{t^2}{n-1}. 
\en
A similar inequality is also true for the Orlicz $\psi_2$-norm of
$f_\theta(t) - f(t)$ generated by the Young function
$\psi_2(r) = e^{r^2} - 1$.

As it turns out, this rate of concentration may be improved
under a second order correlation condition (1.3) at least for 
values of $t$ which are not too large, by involving 
the characteristic $\Lambda = \Lambda(X)$.
In the isotropic case, this condition is described as the relation
\be
\E\,\Big|\sum_{j,k=1}^n z_{jk}\, (X_j X_k - \delta_{jk})\Big|^2
 \, \leq \, \Lambda \sum_{j,k=1}^n |z_{jk}|^2, \quad z_{jk} \in \C.
\en
Here, $\Lambda$ is necessarily bounded away from zero. Indeed, (5.3) includes 
$\E\, X_j^2 X_k^2 - \delta_{jk} \leq \Lambda$ as partial cases. Summing this 
over all $j,k = 1,\dots,n$ leads to $\E\,|X|^4 - n \leq n^2 \Lambda$. But 
$\E\,|X|^4 \geq (\E\,|X|^2)^2 = n^2$ implying that 
$$
\Lambda \geq \frac{n-1}{n} \geq \frac{1}{2}.
$$ 

As was proved  in \cite{B-C-G5} on the basis of Proposition 4.1, if the 
distribution of $X$ is isotropic and symmetric about the origin, the characteristic 
functions $f_\theta(t)$ satisfy in the interval $|t| \leq An^{1/5}$
\be
c\,\E_\theta\, |f_\theta(t) - f(t)|^2\,\leq\,\frac{\Lambda t^4}{n^2}, 
\en
where the constant $c > 0$ depends on the parameter $A \geq 1$ 
only. Moreover,
\be
\E_\theta\, \exp\Big\{\frac{cn}{\Lambda t^2}\, |f_\theta(t) - f(t)|\Big\} 
\leq 2.
\en
Note that, in the symmetric case, the functions 
$\theta \rightarrow f_\theta(t)$ are even, so, all $u_t$ 
have zero linear parts when projecting them onto the subspace 
$H$ of all linear functions in $L^2(\R^n,\s_{n-1})$.

To drop the symmetry assumption, consider an orthogonal decomposition
\be
u_t(\theta) = f(t) + l_t(\theta) + v_t(\theta),
\en
where 
$$
l_t(\theta) = c_1(t)\, \theta_1 + \dots + c_n(t)\, \theta_n, \qquad 
\theta = (\theta_1,\dots,\theta_n) \in \R^n,
$$
is the orthogonal projection of $u_t - f(t)$ onto $H$ (the linear part) and
$v_t(\theta) = u_t(\theta) - f(t) - l_t(\theta)$ is the non-linear part
of $u_t$. By the orthogonality,
\be
\E_\theta\,|f_\theta(t) - f(t)|^2
 \, = \, \E_\theta\, |l_t(\theta)|^2 + \E_\theta\, |v_t(\theta)|^2.
\en
With these notations, the bounds (5.4)-(5.5) should be properly
modified.

\vskip5mm
{\bf Proposition 5.1.} {\sl Given an isotropic random vector $X$ 
in $\R^n$, in the interval $|t| \leq An^{1/5}$,
\be
c\,\E_\theta\, |f_\theta(t) - f(t)|^2\,\leq\, 
\|l_t\|_{L^2}^2 + \frac{\Lambda t^4}{n^2}
\en
with some constant $c > 0$ depending on the parameter $A \geq 1$.
Here $l_t$ is the linear part of $f_\theta(t)$ in $L^2(\R^n,\s_{n-1})$
from the orthogonal decomposition $(5.6)$. Moreover, if 
$|t| \leq An^{1/6}$, then
\be
c\,\|f_\theta(t) - f(t)\|_{\psi_1} \, \leq \, 
\|l_t\|_{L^2} + \frac{\Lambda t^2}{n}.
\en
}

\vskip2mm
If the distribution of $X$ is symmetric about the origin, then
$l_t(\theta) = 0$, and we return in (5.8)-(5.9) to (5.4)-(5.5). 
The linear part $l_t$ is also vanishing, when $X$ has mean zero 
and a constant Euclidean norm, i.e. when $|X| = \sqrt{n}$ a.s. 
(this will be clarified in the next section). 

\vskip5mm
{\bf Proof.} To employ Propositions 4.1-4.2, we need to choose 
a suitable value $a \in \mathbb C$ and estimate the operator norm 
$\|\nabla^2 u_t - a\,\I_n\|$ and the Hilbert-Schmidt norm 
$\|\nabla^2 u_t - a\,\I_n\|_{\rm HS}$. First note that, 
by differentiation of (5.1), for any fixed $t \in \R$,
$$
\big[\nabla^2 u_t(\theta)\big]_{jk} =
\frac{\partial^2}{\partial \theta_j \partial \theta_k}\,f_\theta(t) =
-t^2\, \E\,X_j X_k\,e^{it\left<X,\theta\right>}.
$$
Hence, a good choice is $a = -t^2 f(t)$ in order to balance 
the diagonal elements in the matrix of second derivatives of $u_t$. 
For any vector $w \in {\mathbb C}^n$, using 
the canonical inner product in the complex $n$-space, we have
$$
\left<\nabla^2 u_t(\theta)\,w,w\right> = -t^2\,
\E\,|\left<X,w\right>|^2\,e^{it\left<X,\theta\right>}.
$$
Hence, by the isotropy assumption,
$$
\big|\left<(\nabla^2 u_t(\theta) - a\,\I_n)\,w,w\right>\big| \leq t^2\,
\E\,|\left<X,w\right>|^2 + |a|\,|w|^2 \leq 2t^2, \qquad |w|=1.
$$
In terms of the norm defined as in (4.8), this bound insures that 
\be
\|\nabla^2 u_t(\theta) - a\,\I_n\| \leq 2t^2.
\en

In addition, putting $a(\theta) = -t^2 f_\theta(t)$, we have
\bee
\big\|\nabla^2 u_t(\theta) - a(\theta)\, \I_n\big\|_{\rm HS}^2 
 & = &
\sum_{j,k = 1}^n \left|\nabla^2 u_t(\theta)_{jk} - 
a(\theta)\,\delta_{jk}\right|^2 \\
 & = &
\sup \bigg|\sum_{j,k=1}^n z_{jk}\,
\big(\nabla^2 u_t(\theta)_{jk} - a(\theta)\,\delta_{jk}\big)\bigg|^2 \\ 
 & \leq &
t^4 \,\sup\,  
\E\,\bigg|\sum_{j,k=1}^n z_{jk}\, (X_j X_k - \delta_{jk})\bigg|^2,
\ene
where the supremum is running over all complex numbers $z_{jk}$ such that 
$\sum_{j,k=1}^n |z_{jk}|^2 = 1$. But, under this constraint, due to 
the second order correlation condition, the last
expectation is bounded by $\Lambda$. Since $u_t$ and $v_t$ 
have equal Hessians, we conclude that
\be
\big\|\nabla^2 v_t(\theta) - a(\theta)\, \I_n\big\|_{\rm HS}^2 \, \leq \,
\Lambda t^4
\en
for all $\theta$. On the other hand, by (5.2),
\be
\E_\theta\, \big\|(a(\theta) - a)\, \I_n\big\|_{\rm HS}^2 
 \, = \, 
n t^4\,\E_\theta\,|f_\theta(t) - f(t)|^2 \, \leq \, 2t^6.
\en
The two last bounds give
$$
\E_\theta\, \big\|\nabla^2 v_t(\theta) - a\,\I_n\big\|_{\rm HS}^2 
 \, \leq \, 2\Lambda t^4 + 4t^6,
$$
which, by Proposition 4.1, yields
$$
\E_\theta\, |v_t(\theta)|^2 \, \leq \, 
\frac{5}{(n-1)^2}\,(2\Lambda t^4 + 4 t^6).
$$

One can sharpen this bound for the range 
$|t| \leq An^{1/5}$. Applying it in (5.7),~we~get
$$
\E_\theta\,|f_\theta(t) - f(t)|^2 \, \leq \, 
\E_\theta\, |l_t(\theta)|^2 + \frac{5}{(n-1)^2}\,(2\Lambda t^4 + 4 t^6),
$$
which, according to the identity in (5.12), gives
$$
\E_\theta \|(a(\theta) - a)\, \I_n\|_{\rm HS}^2 \, \leq \,
nt^4\,\E_\theta\, |l_t(\theta)|^2 + 
\frac{5n}{(n-1)^2}\,(2\Lambda t^8 + 4 t^{10}).
$$
Combining this with (5.11), we get
$$
\E_\theta \|\nabla^2 v_t(\theta) - a\, \I_n\|_{\rm HS}^2 \, \leq \, 
2nt^4\,\E_\theta\, |l_t(\theta)|^2 + 2\Lambda t^4 +
\frac{10\,n}{(n-1)^2}\,(2\Lambda t^8 + 4 t^{10}).
$$
Hence, by Proposition 4.1 once more,
$$
\E_\theta\, |v_t(\theta)|^2 \, \leq \, 
\frac{10\,n t^4}{(n-1)^2}\, \E_\theta\, |l_t(\theta)|^2 + 
\frac{10}{(n-1)^2}\,\Lambda t^4 + 
\frac{50\,n}{(n-1)^4}\,(\Lambda t^8 + 2 t^{10}),
$$
so that, by (5.7),
\bee
\E_\theta\,|f_\theta(t) - f(t)|^2 
 & \leq &
\Big(1 + \frac{10\,nt^4}{(n-1)^2}\Big)\, \E_\theta\, |l_t(\theta)|^2 \\
 & & + \
\frac{10}{(n-1)^2}\,\Lambda t^4 + 
\frac{50\,n}{(n-1)^4}\,(\Lambda t^8 + 2 t^{10}).
\ene

According to the identity in (15.12), this gives
\bee
\E_\theta \|(a(\theta) - a)\, \I_n\|_{\rm HS}^2 
 & \leq &
nt^4\,\Big(1 + \frac{10\,n t^4}{(n-1)^2}\Big)\, 
\E_\theta\, |l_t(\theta)|^2 \\
 & & + \
\frac{10\,n}{(n-1)^2}\,\Lambda t^8 + 
\frac{50\,n^2}{(n-1)^4}\,(\Lambda t^{12} + 2 t^{14}).
\ene
One can combine this with (5.11) to obtain that
\bee
\E_\theta\, \|\nabla^2 v_t(\theta) - a\, \I_n\|_{\rm HS}^2 
 & \leq & 
2nt^4\,\Big(1 + \frac{10\,n t^4}{(n-1)^2}\Big)\, 
\E_\theta\, |l_t(\theta)|^2 \\
 & & + \  
2\Lambda t^4 + \frac{20\,n}{(n-1)^2}\,\Lambda t^8 + 
\frac{100\,n^2}{(n-1)^4}\,(\Lambda t^{12} + 2 t^{14}).
\ene

Now, if $|t| \leq An^{1/5}$, the coefficient in front of 
$\E_\theta\, |l_t(\theta)|^2$ does not exceed a multiple of 
$nt^4$. Similarly, in this region the last three terms can be 
bounded by $\Lambda t^4$ up to a numerical factor (since 
$\Lambda \geq \frac{1}{2}$). Hence the above bound is simplified to
\be
c\,\E_\theta\, \|\nabla^2 v_t(\theta) - a\, \I_n\|_{\rm HS}^2 
 \, \leq \, nt^4\,\|l_t\|_{L^2}^2 + \Lambda t^4
\en
with some constant $c$ depending $A$.
Since $nt^4 < A^4 n^2$, by Proposition 4.1, we get
$$
c\,\E_\theta\, |v_t(\theta)|^2 \, \leq \, 
\E_\theta\, |l_t(\theta)|^2 + \frac{5}{(n-1)^2}\, \Lambda t^4.
$$
In view of (5.7), this proves the inequality (5.8).

To get a bound for the $\psi_1$-norm, note that, by (5.10), 
the conditions of Proposition 4.2 (in its second part) are fulfilled 
with $-\frac{1}{2}\,f(t)$ in place of $a$ for the function
$$
u(\theta) = \frac{1}{2t^2}\,(f_\theta(t) - f(t)), \qquad 
\theta \in \R^n, \ \ t \neq 0.
$$
Since (5.13) holds for $u_t$ as well (provided that 
$|t| \leq An^{1/5}$), this inequality may be rewritten as
$$
c\,\E_\theta\, 
\Big\|\nabla^2 u(\theta) + \frac{1}{2}\,f(t)\, \I_n\Big\|_{\rm HS}^2 
 \, \leq \, n\,\|l_t\|_{L^2}^2 + \Lambda.
$$
The linear part of $u$ is given by $l_t/(2t^2)$.
Hence, the inequality (4.6) of Proposition 4.2 yields
$$
c\, \Big\|\frac{1}{2t^2}\,(f_\theta(t) - f(t))\Big\|_{\psi_1} \, \leq \, 
\frac{1 + \Lambda}{n} + \|l_t\|_{L^2}^2 + \frac{1}{2t^2}\,\|l_t\|_{L^2}.
$$
Using once more $\Lambda \geq \frac{1}{2}$, the above is simplified
to 
\be
c\,\|f_\theta(t) - f(t)\|_{\psi_1} \, \leq \, \frac{\Lambda t^2}{n} +
\|l_t\|_{L^2} + \|l_t\|_{L^2}^2\, t^2.
\en

Here, the last term on the right-hand side is dominated by the second
last term in the smaller interval $|t| \leq An^{1/6}$. Indeed,
according to the concentration inequality (5.2),
$$
\|l_t\|_{L^2}\,t^2 \, \leq \, \|f_\theta(t) - f(t)\|_{L^2}\,t^2 \, \leq \,
\frac{|t|^3}{\sqrt{n-1}} \, < \, 2A^3.
$$
As a result, (5.14) leads to the required form (5.9).
\qed

\vskip5mm
{\bf Remark.}
Continuing the iteration process in the proof of Proposition 5.2, 
one may state (5.8) in the intervals $|t| \leq n^\alpha$ with 
any fixed $\alpha < \frac{1}{4}$.


\vskip7mm
\noindent
\section{{\bf Linear Part of Characteristic Functions}}
\setcounter{equation}{0}

\vskip2mm
\noindent
In order to make the bounds (5.8)-(5.9) effective, we need to properly 
estimate the $L^2$-norm of the linear part $l_t(\theta)$ of $f_\theta(t)$ 
in $L^2(\R^n,\s_{n-1})$. According to (4.10), it is described as
\be
I(t) \, = \, \|l_t\|_{L^2}^2 \, = \,
n\,\E_\theta\,\E_{\theta'} \left<\theta,\theta'\right> 
f_\theta(t) \bar f_{\theta'}(t).
\en
Let us find an asymptotically explicit expression for this function.

\vskip5mm
{\bf Proposition 6.1.} {\sl Let $X$ be a random vector in $\R^n$ 
such that $\E\,|X|^2 = n$. For any $t \in \R$, the characteristic function
$f_\theta(t) = \E\,e^{it\left<X,\theta\right>}$ as a function of 
$\theta$ on the sphere has a linear part, whose squared $L^2(\s_{n-1})$-norm 
may be represented as
\be
I(t) = \frac{t^2}{n}\, \E \left<X,Y\right>
\Big(1 - \frac{(U^2 + V^2)\, t^4 - 8R^2 t^2}{4n}\Big)\,
e^{-R^2 t^2} + O(t^2 n^{-2}),
\en
where $Y$ is an independent copy of $X$ and
$R^2 = \frac{|X|^2 + |Y|^2}{2n}$, $U = \frac{|X|^2}{n}$, 
$V = \frac{|Y|^2}{n}$. 
The remainder term may be improved to $O(t^2 n^{-5/2})$,
if $X$ is isotropic.
}

\vskip5mm
{\bf Proof.} Using an independent copy $Y$ of $X$, one may
rewrite (6.1) equivalently as
$$
I(t) \, = \, n\,\sum_{k=1}^n |\,\E_\theta\,\theta_k f_\theta(t)|^2 \, = \, 
n \sum_{k=1}^n\, \E\ \E_\theta\,\E_{\theta'} \Big[\theta_k \theta_k'\,
e^{it\left<X,\theta\right> - it\left<Y,\theta'\right>}\Big].
$$
To compute the inner expectations, introduce the function
$$
K_n(t) = J_n\big(\sqrt{tn}\big), \qquad t \geq 0,
$$
where, as before, $J_n$ denotes the characteristic function of the 
first coordinate of a point on the unit sphere $\S^{n-1}$ under the
normalized Lebesgue measure $\s_{n-1}$. By the definition,
$$
\E_\theta\, e^{i\left<v,\theta\right>} = J_n(|v|) =
K_n\Big(\frac{|v|^2}{n}\Big), \qquad v=(v_1,\dots,v_n) 
\in \R^n.
$$
Differentiating this equality with respect to the variable $v_k$, 
we obtain that
$$
i\,\E_\theta\, \theta_k e^{i\left<v,\theta\right>} \, = \,
\frac{2v_k}{n}\,K_n'\Big(\frac{|v|^2}{n}\Big).
$$
Let us multiply this by a similar equality
$$
-i\,\E_\theta\, \theta_k e^{-i\left<w,\theta\right>} =
\frac{2w_k}{n}\,K_n'\Big(\frac{|w|^2}{n}\Big),
$$
to get that, for all $v,w \in \R^n$,
$$
\E_\theta\,\E_{\theta'} \Big[\theta_k \theta_k'\,
e^{i\left<v,\theta\right> - i\left<w,\theta'\right>}\Big] = 
\frac{4v_k w_k}{n^2}\, 
K_n'\Big(\frac{|v|^2}{n}\Big)\,K_n'\Big(\frac{|w|^2}{n}\Big).
$$
Hence, summing over all $k \leq n$, we get
$$
\sum_{k=1}^n \E_\theta\,\E_{\theta'} \Big[\theta_k \theta_k'\,
e^{i\left<v,\theta\right> - i\left<w,\theta'\right>}\Big] = 
\frac{4\left<v,w\right>}{n^2}\,
K_n'\Big(\frac{|v|^2}{n}\Big)\,K_n'\Big(\frac{|w|^2}{n}\Big).
$$
It remains to make the substitution $v = tX$, $w = tY$ and to 
take the expectation over $(X,Y)$. Then we arrive at the following 
expression 
\be
I(t) = \frac{4t^2}{n}\, \E \left<X,Y\right>
K_n'\Big(\frac{t^2 |X|^2}{n}\Big) K_n'\Big(\frac{t^2 |Y|^2}{n}\Big).
\en

In particular, if $|X| = \sqrt{n}$ a.s., then 
$$
I(t) = \frac{4t^2}{n}\,K_n'^2(t^2) \, \E \left<X,Y\right>,
$$
which is vanishing, as soon as $X$ has mean zero. 
In fact, the property $I(t) = 0$ remains valid for more
general random vectors. In particular, this is the case, where 
the conditional distribution of $X$ given that $|X| = r$ has mean 
zero for any $r>0$.

Now, let us derive an asymptotic formula for the function 
$K_n$ and its derivative. We know from Corollary 3.2 that
$$
\frac{d}{dt}\, J_n(t\sqrt{n}) = 
-t\Big(1 - \frac{t^4 - 4t^2}{4n}\Big)\,e^{-t^2/2} + 
O\big(n^{-2} \min(1,|t|^3)\big).
$$
Since $J_n(t\sqrt{n}) = K_n(t^2)$, after differentiation we find that
$$
2t K_n'(t^2) = \frac{d}{dt}\, K_n(t^2) = 
-t\Big(1 - \frac{t^4 - 4t^2}{4n}\Big)\,e^{-t^2/2} + 
O\big(n^{-2} \min(1,|t|^3)\big).
$$
Changing the variable, we arrive at
$$
K_n'(t) = -\frac{1}{2}\,
\Big(1 - \frac{t^2 - 4t}{4n}\Big)\,e^{-t/2} + O\big(n^{-2} \min(1,t)\big),
\qquad t \geq 0.
$$

From this,
$$
K_n'(t)K_n'(s) \, = \, \frac{1}{4}\,
\Big(1 - \frac{(t^2 + s^2) - 4(t+s)}{4n}\Big)\,e^{-(t+s)/2} + 
O\big(n^{-2}\big)
$$
uniformly over all $t,s \geq 0$, so,
\bee
4 K_n'\Big(\frac{t^2 |X|^2}{n}\Big) K_n'\Big(\frac{t^2 |Y|^2}{n}\Big)
 & = &
\bigg(\!1 - \frac{t^4\,(\frac{|X|^4}{n^2} + \frac{|Y|^4}{n^2}) 
- 4t^2\,(\frac{|X|^2}{n} + \frac{|Y|^2}{n})}{4n}\bigg)\,
e^{-\frac{t^2 (|X|^2 + |Y|^2)}{2n}} + \ep  \\
 & = &
\Big(1 - \frac{(U^2 + V^2)\, t^4 - 8R^2 t^2}{4n}\Big)\,
e^{-R^2 t^2} + \ep
\ene
with a remainder term satisfying $|\ep| \leq \frac{c}{n^2}$
up to some absolute constant $c$. The latter yields
$$
\frac{t^2}{n}\,\E\,|\left<X,Y\right>|\,|\ep| \, \leq \,
\frac{ct^2}{n^2}\,\E\, \frac{|X|^2 + |Y|^2}{2n} \, = \, \frac{ct^2}{n^2},
$$
assuming that $\E\,|X|^2 = n$. Hence, recalling (6.3), we obtain (6.2).

In the isotropic case, we have $\E\,|\left<X,Y\right>| \leq \sqrt{n}$, which
leads to the corresponding improvement of the remainder term.
\qed

 
\vskip7mm 
\section{{\bf Berry-Esseen Bounds}} 
\setcounter{equation}{0} 
 
\vskip2mm 
\noindent 
The Kolmogorov distances between the distribution functions $F_\theta$ 
of the weighted sums $S_\theta = \left<X,\theta\right>$ and the standard 
normal distribution function $\Phi$ can be explored by means of the 
Berry-Esseen-type bounds. They involve the characteristic functions
\be
f_\theta(t) = \E\,e^{itS_\theta} =  
\int_{-\infty}^\infty e^{itx}\, dF_\theta(x), \qquad
f(t) = \E_\theta f_\theta(t)= \int_{-\infty}^\infty e^{itx}\, dF(x) 
\en
associated to $F_\theta(x)$ and the average distribution function 
$F(x) = \E_\theta F(x)$. Using the $\Lambda$-functional, 
let us state a few preliminary relations.

\vskip5mm
{\bf Lemma 7.1.} {\sl Given a random vector $X$ in $\R^n$ such that 
$\E\, |X|^2 = n$, we have, for all $T \geq T_0 \geq 1$ and 
$\theta \in \S^{n-1}$, 
\begin{eqnarray}
c\,\rho(F_\theta,\Phi) 
 & \leq &
\int_0^{T_0} \frac{|f_\theta(t) - f(t)|}{t}\,dt \nonumber \\
 & & + \
\int_{T_0}^T \frac{|f_\theta(t)|}{t}\,dt  + \frac{\Lambda}{n}
\,\Big(1 + \log\frac{T}{T_0}\Big) + \frac{1}{T} + e^{-T_0^2/4}.
\end{eqnarray}
}

\vskip2mm
The idea to involve two parameters $T$ and $T_0$ stems upon the observation
that the first integrand in (7.2) is small on a relatively moderate sized interval
$[0,T_0]$ only, due to the concentration property of $f_\theta(t)$ about $f(t)$ 
as a function of $\theta$ (as discussed in Section 5). On the other hand,
for $T_0 \leq t \leq T$ with a sufficiently large $T$, one may hope that 
both $f_\theta(t)$ and $f(t)$ will be just small in absolute value 
(in analogy with the case of independent components). 

\vskip5mm
{\bf Proof.} One can apply a general Berry-Esseen-type bound
$$
c\,\rho(U,V) \, \leq \,
\int_0^T \frac{|\hat U(t) - \hat V(t)|}{t}\,dt + \frac{1}{T}  
\int_0^T |\hat V(t)|\,dt \qquad (T>0),
$$
where $U$ and $V$ are arbitrary distribution functions
with characteristic functions $\hat U$ and $\hat V$, respectively 
(cf. e.g. \cite{B4}, \cite{P1}, \cite{P2}). In particular, for all
$\theta \in \S^{n-1}$,
$$
c\,\rho(F_\theta,F) \, \leq \,
\int_0^T \frac{|f_\theta(t) - f(t)|}{t}\,dt + \frac{1}{T}  
\int_0^T |f(t)|\,dt.
$$
Splitting the integration in the first integral to the subintervals
$[0,T_0]$ and $[T_0,T]$, $T \geq T_0 > 0$, we then have
\begin{eqnarray}
c\,\rho(F_\theta,F) 
 & \leq &
\int_0^{T_0} \frac{|f_\theta(t) - f(t)|}{t}\,dt \nonumber \\ 
 & & + \
\int_{T_0}^T \frac{|f_\theta(t)|}{t}\,dt +
\int_{T_0}^T \frac{|f(t)|}{t}\,dt + \frac{1}{T} \int_0^T |f(t)|\,dt.
\end{eqnarray}

The decay of the characteristic function $f(t)$ for large $t$ 
can be controlled in terms of the variance-type functional 
$\sigma_4^2 = \frac{1}{n}\,\Var(|X|^2)$, which in turn satisfies 
$\sigma_4^2 \leq \Lambda$ according to the inequality (1.3) applied with 
coefficients $a_{ij} = 1$. Namely, write the definition (7.1) as
$$
f(t) = \E\,J_n(t|X|), \qquad t \in \R.
$$
Here, one may split the expectation into the event 
$A = \{|X|^2 \leq \frac{1}{2}\,n\}$ and its complement $B$.
By the upper bound (3.1),
$$
\E\,|J_n(t|X|)|\,1_B \, \leq \, 
\E\,\Big(5\,e^{-|X|^2/2} + 4\,e^{-n/12}\Big)\,1_B \, \leq \, 
5\,e^{-t^2/4} + 4\,e^{-n/12}.
$$
On the other hand, by Chebyshev's inequality,
\be
\P(A) = \P\Big\{n - |X|^2 \geq \frac{1}{2}\,n\Big\} 
\, \leq \, \frac{\Var(|X|^2)}{(\frac{1}{2}\,n)^2} \, = \,
\frac{4\sigma_4^2}{n} \, \leq \, \frac{4\Lambda}{n}.
\en
Since $|J_n(s)| \leq 1$ for all $s \in \R$, we get
$$
\E\,|J_n(t|X|)|\,1_A \, \leq \, \frac{4\Lambda}{n},
$$
thus implying that $c\,|f(t)| \leq e^{-t^2/4} + \frac{\Lambda}{n}$ 
for all $t \in \R$, and therefore
\be
\frac{c}{T} \int_0^T |f(t)|\,dt \, \leq \,  
\frac{\Lambda}{n} + \frac{1}{T}.
\en

If $T_0 \geq 1$, then also
\be
c \int_{T_0}^T \frac{|f(t)|}{t}\,dt \, \leq \, e^{-T_0^2/4} + 
\frac{\Lambda}{n}\,\log(T/T_0).
\en
Using these bounds in the inequality (7.3), it is simplified to
\bee
c\,\rho(F_\theta,F) 
 & \leq &
\int_0^{T_0} \frac{|f_\theta(t) - f(t)|}{t}\,dt \\
 & & + \
\int_{T_0}^T \frac{|f_\theta(t)|}{t}\,dt  + \frac{\Lambda}{n}
\,\Big(1 + \log\frac{T}{T_0}\Big) + \frac{1}{T} + e^{-T_0^2/4}.
\ene

The variance functional may also be used to quantify closeness of $F$
to the standard normal distribution function via the inequality
(cf. \cite{B-C-G3})
$$
c\,\rho(F,\Phi) \leq \frac{1}{n}\,(1 + \sigma_4^2), 
$$
Since $\sigma_4^2 \leq \Lambda$, (7.2) immediately follows
in view of the triangle inequality for the Kolmogorov metric.
\qed

\vskip5mm
Lemma 7.1 may be used to derive the following upper bound on average
which represents a generalization of the inequality (1.4).

\vskip5mm
{\bf Lemma 7.2.} {\sl Given an isotropic random vector $X$ in $\R^n$, 
with $T_0 = 4\sqrt{\log n}$ we have
\be
c\,\E_\theta\, \rho(F_\theta,\Phi) \leq \frac{\log n}{n}\, \Lambda +
\int_0^{T_0} \frac{\sqrt{I(t)}}{t}\,dt,
\en
where $I(t)$ denotes the squared $L^2$-norm of the linear part of 
$f_\theta(t)$ in $L^2(\s_{n-1})$.
}

\vskip5mm
{\bf Proof.}
When bounding $\rho(F_\theta,\Phi)$ on average with respect to
$\s_{n-1}$, the inequality (7.6) is actually not needed. Using Jensen's 
inequality $|f(t)| \leq \E_\theta\, |f_\theta(t)|$,
from (7.3) and (7.5) we obtain that, for all $T \geq T_0 \geq 1$,
\be
c\,\E_\theta\, \rho(F_\theta,F) \, \leq \,
\int_0^{T_0} \frac{\E_\theta\,|f_\theta(t) - f(t)|}{t}\,dt +
\int_{T_0}^T \frac{\E_\theta\,|f_\theta(t)|}{t}\,dt + 
\frac{1}{T}  + \frac{\Lambda}{n}.
\en

Now, as was shown in \cite{B-C-G4} (Lemma 5.2 specialized 
to the parameter $p = 2$), for all $t \in \R$,
\be
c\, \E_\theta\,|f_\theta(t)| \, \leq \, \frac{m_4^2 + \sigma_4^2}{n} +
e^{-t^2/16}, \quad
m_4 = \frac{1}{\sqrt{n}}\,\big(\E \left<X,Y\right>^4\big)^{1/4},
\en
where $Y$ is an independent copy of $X$. Using a simple relation $m_4 \leq M_4^2$ 
(Corollary 2.3 in \cite{B-C-G4}), one may also involve the functional
$$
M_4 = \sup_{\theta \in \S^{n-1}}\,(\E \left<X,\theta\right>^4)^{1/4}.
$$
It may be bounded in terms of $\Lambda$ as well as $\sigma_4^2$.
Indeed, applying (1.3) with $a_{ij} = \theta_i \theta_j$, we get
$$
\Var(\left<X,\theta\right>^2) \leq \Lambda, \quad \theta \in S^{n-1},
$$
which implies $M_4^4 \leq 1 + \Lambda \leq 3\Lambda$ in the isotropic case. 
This allows us to replace (7.9) with
$$
c\, \E_\theta\,|f_\theta(t)| \, \leq \, \frac{\Lambda}{n} +
e^{-t^2/16}.
$$
Applying the latter in (7.8), thus inequality is simplified to
\be
c\,\E\,\rho(F_\theta,\Phi) \, \leq \,
\int_0^{T_0} \frac{\E\,|f_\theta(t) - f(t)|}{t}\,dt + \frac{\Lambda}{n}
\,\Big(1 + \log\frac{T}{T_0}\Big) + \frac{1}{T} + e^{-T_0^2/16}.
\en

Here, the integral can be bounded by virtue of the $L^2$-bound (5.8) which
yields
$$
c\,\E_\theta\, |f_\theta(t) - f(t)| \, \leq \,
\sqrt{I(t)} + \frac{t^2}{n}\sqrt{\Lambda}
$$
for $|t| \leq An^{1/5}$ with a prescribed constant $A \geq 1$. This gives
$$
c \int_0^{T_0} \frac{\E\,|f_\theta(t) - f(t)|}{t}\,dt \, \leq \,
\int_0^{T_0} \frac{\sqrt{I(t)}}{t}\,dt + \frac{T_0^2}{2n}\sqrt{\Lambda},
$$
as long as $T_0 \leq An^{1/5}$. Applying this in (7.10), we arrive at
$$
c\,\E_\theta\,\rho(F_\theta,\Phi) 
 \, \leq \,
\int_0^{T_0} \frac{\sqrt{I(t)}}{t}\,dt +
\frac{\Lambda}{n}\,\Big(1 + \log\frac{T}{T_0}\Big) + 
\frac{1}{T} + \frac{T_0^2}{n}\sqrt{\Lambda} + e^{-T_0^2/16}.
$$
Finally, choosing $T = 4n$, $T_0 = 4\sqrt{\log n}$, we obtain (7.7).
\qed


\vskip7mm
\section{{\bf Large Deviations Related to Moderate Sized and Long Intervals}}
\setcounter{equation}{0}

\vskip2mm
\noindent
A similar argument can be used when 
bounding the $\psi_1$-Orlicz norm of $\rho(F_\theta,\Phi)$. 
As a preliminary step, let us start with the first integral in (7.2) 
over the moderate interval. Applying now the inequality (5.9), we have
\bee
c\, \bigg\|\int_0^{T_0}|f_\theta(t) - f(t)|\,\frac{dt}{t}\,\bigg\|_{\psi_1} 
 & \leq & 
c \int_0^{T_0} \|f_\theta(t) - f(t)\|_{\psi_1}\,\frac{dt}{t} \\
 & \leq & 
\int_0^{T_0} \Big(\sqrt{I(t)} + \frac{\Lambda t^2}{n}\Big)\,\frac{dt}{t} \\
 & = &
\frac{\Lambda T_0^2}{2n} + \int_0^{T_0} \frac{\sqrt{I(t)}}{t}\,dt,
\ene
which is used with the same parameter $T_0$ as in Lemma 7.2.
In general, by Markov's inequality, 
$$
\s_{n-1}\big\{|\xi| \geq r\|\xi\|_{\psi_1}\big\} \leq 2\,e^{-r}, \quad r>0.
$$ 
Hence, we get:

\vskip5mm
{\bf Lemma 8.1.} {\sl Let $X$ be an isotropic random 
vector in $\R^n$. For all $r > 0$, with $T_0 = 4\sqrt{\log n}$,
$$
\s_{n-1}\bigg\{c \int_0^{T_0} \frac{|f_\theta(t) - f(t)|}{t}\,dt 
\geq \frac{\Lambda \log n}{n}\,r + r
\int_0^{T_0} \frac{\sqrt{I(t)}}{t}\,dt\bigg\} \, \leq \, 2\,e^{-r}.
$$
}

\vskip2mm
Outside the moderate sized interval, that is, on the long interval 
$[T_0,T]$, both $|f(t)|$ and $|f_\theta(t)|$ are expected to be small 
for most of $\theta$. To study this property, let us consider the growth of
the moments of the integral
\be
L(\theta) = \int_{T_0}^T \frac{|f_\theta(t)|}{t}\,dt.
\en

\vskip5mm
{\bf Lemma 8.2.} {\sl Given a random vector $X$ in $\R^n$, let
$X^{(k)}$, $Y^{(k)}$ $(k = 1,\dots,p)$ be independent copies of $X$.
For the integral in $(8.1)$ with parameters 
$T_0 = 4\sqrt{\log n}$ and $T = T_0 n$, we have
\be
\E_\theta\,L(\theta)^{2p} \, \leq \, (c\log n)^{2p}\, 
\big(p^{2p}\, n^{-2p} + \P(A)\big),
\en
where 
$$
A = \Big\{|\Sigma_p|^2 \leq \frac{np}{2}\Big\}, \quad
\Sigma_p = \sum_{k=1}^p\, (X^{(k)} - Y^{(k)}).
$$
}

\vskip2mm
{\bf Proof.} By H\"older's inequality,
$$
L(\theta)^{2p} \, \leq \, \log^{2p-1}\Big(\frac{T}{T_0}\Big)\, 
\int_{T_0}^T \frac{|f_\theta(t)|^{2p}}{t}\,dt,
$$
so that
$$
\E_\theta\,L(\theta)^{2p} \, \leq \, 
\log^{2p-1}\Big(\frac{T}{T_0}\Big)\, 
\int_{T_0}^T \frac{\E_\theta\,|f_\theta(t)|^{2p}}{t}\,dt.
$$
Since $|f_\theta(t)|^{2p} = \E\, e^{it\left<\Sigma_p,\theta\right>}$,
we may write
$$
\E_\theta\,|f_\theta(t)|^{2p} \, = \, \E J_n(t\,|\Sigma_p|).
$$
Thus,
$$
\E_\theta\,L(\theta)^{2p} \, \leq \, 
\log^{2p-1}\Big(\frac{T}{T_0}\Big)\, 
\int_{T_0}^T \E J_n(t\,|\Sigma_p|)\,\frac{dt}{t}.
$$

Next, we split the expectation to the events $A$ and its complement 
$B = \big\{|\Sigma_p|^2 > \frac{np}{2}\big\}$. 
Applying the upper bound (3.1), we get
\bee
\int_{T_0}^T \E J_n(t\,|\Sigma_p|)\,1_B\,\frac{dt}{t} 
 & \leq & 
\int_{T_0}^T \frac{5\,e^{-pt^2/4} + 4\,e^{-n/12}}{t}\,dt \\
 & \leq &
(5\,e^{-pT_0^2/2} + 4\,e^{-n/12})\log\Big(\frac{T}{T_0}\Big),
\ene
while
$$
\int_{T_0}^T \E J_n(t|\Sigma_p|)\,1_A\,\frac{dt}{t} 
 \, \leq \, \P(A) \log\Big(\frac{T}{T_0}\Big)
$$
(since $|J(s)| \leq 1$ for all $s \in \R$). Hence,
$$
\E_\theta\,L(\theta)^{2p} \, \leq \, c \log^{2p}\Big(\frac{T}{T_0}\Big)\, 
\Big(e^{-pT_0^2/2} + e^{-n/12} + \P(A)\Big).
$$
For the choice $T_0 = 4\sqrt{\log n}$, $T = T_0 n$, this leads to
$$
\E_\theta\,L(\theta)^{2p} \, \leq \, c\,(\log n)^{2p}\, 
\big(n^{-8p} + e^{-n/12} + \P(A)\big).
$$
Using the inequality 
$x^{2p}\, e^{-x} \leq p^{2p}$ ($x \geq 0$), 
we have $e^{-n/12} \leq (12\, p)^{2p}\,n^{-2p}$, 
and the above bound is simplified to (8.2).
\qed


\vskip7mm
\section{{\bf Concentration in Presence of Poincar\'e-type Inequalities}}
\setcounter{equation}{0}

\vskip2mm
\noindent
In order to simplify the bounds in Lemma 7.2 and Lemmas 8.1-8.2,
we need more information about the distribution of $X$, which would 
allow us to say more on the involved function $I(t)$ and
the probability of the even $A$ as in Lemma 8.2.
To this aim, our starting hypothesis will be described by 
Poincar\'e-type inequalities.

Let us first recall several results about concentration,
assuming that the random vector $X = (X_1,\dots,X_n)$ in $\R^n$
admits the Poincar\'e-type inequality 
\be
\lambda_1\, \Var(u(X)) \, \leq \, \E\, |\nabla u(X)|^2
\en
for all smooth functions $u$ on $\R^n$ with a positive constant 
$\lambda_1$. As was discovered by Gromov and Milman \cite{G-M} and 
by Borovkov and Utev \cite{B-U}, deviations of random variables $u(X)$ 
from their means are subexponential, as long as $u$ is a Lipschitz 
function on $\R^n$ (cf. also \cite{A-S}, \cite{L2}). In a somewhat 
optimal way, worst possible deviations of $u(X)$ are described 
in the following assertion proved in \cite{B1}.

\vskip5mm
{\bf Proposition 9.1.} {\sl If the function $u:\R^n \rightarrow \R$ 
has a Lipschitz semi-norm $\|u\|_{\rm Lip} \leq 1$, then, for any 
$r \geq 0$,
\be
\P\big\{u(X) - \E\, u(X) \geq r\big\} \, \leq \, 
3\,e^{-2\sqrt{\lambda_1} r}.
\en
}

\vskip2mm
Using a smoothing argument, the inequality (9.2) may be extended to all
locally Lipschitz functions, in which case the modulus of the gradient 
is understood as a Borel measurable function
$$
|\nabla u(x)| \, = \, \limsup_{y \rightarrow x} \frac{|u(x) - u(y)|}{|x-y|},
\quad x \in \R.
$$
In terms of partial derivatives, it leads to the usual expression
$\big(\sum_{k=1}^n (\partial_{x_k} u(x))^2\big)^{1/2}$ assuming that
$u$ is differentiable at the point $x$.

If the function $u$ is not Lipschitz (for example, a polynomial), 
the bound (9.2) is no longer true, and a more general variant 
of Proposition 9.1 is needed, which would allow us to control probabilities 
of large deviations. To this aim, proper bounds on the $L^p$-norms of $u$ in 
terms of the $L^p$-norms of the modulus of the gradient are useful.

\vskip5mm
{\bf Proposition 9.2.} {\sl Given a locally Lipschitz function $u$ on $\R^n$, 
suppose that the moment $\E\, |\nabla u(X)|^p$ is finite for $p \geq 2$. 
Then, $u(X)$ has finite absolute moments up to order $p$, and
\be
\E\, |u(X) - \E\,u(X)|^p \, \leq \,
\Big(\frac{p}{\sqrt{2\lambda_1}}\Big)^p\,\E\, |\nabla u(X)|^p.
\en
}

\vskip2mm
{\bf Proof.} Let us include a simple argument, assuming that the 
function $u$ is $C^1$-smooth. By the subadditivity property of the 
variance functional (cf. \cite{L1}), the Poincar\'e-type inequality (9.1)
for the distribution $\mu$ of $X$ on $\R^n$ is extended 
to the same relation on $\R^n \times \R^n$
\be
\lambda_1 \Var_{\mu \otimes \mu}(f) \leq 
\int\!\!\!\int |\nabla f(x,y)|^2\,d\mu(x)d\mu(y)
\en
with respect to the product measure $\mu^2 = \mu \otimes \mu$. Here, for any 
$C^1$-smooth function $f = f(x,y)$, the modulus of the gradient is given by
$$
|\nabla f(x,y)|^2 = |\nabla_x f(x,y)|^2 + |\nabla_y f(x,y)|^2.
$$
Let us apply this $2n$-dimensional Poincar\'e-type inequality to the function 
$$
f(x,y) = |u(x) - u(y)|^{\frac{p}{2}}\ {\rm sign}(u(x) - u(y)),
$$
which is $C^1$-smooth in case $p > 2$. Its modulus of the gradient
is given by
$$
|\nabla f(x,y)| \, = \, \frac{p}{2}\ |u(x) - u(y)|^{\frac{p}{2} - 1}
\sqrt{|\nabla u(x)|^2 + |\nabla u(y)|^2}.
$$
Since $f$ has a symmetric distribution under $\mu^2$, applying (9.4) 
together with H\"older's inequality, we conclude that
\bee
\lambda_1 \int\!\!\!\int |u(x) - u(y)|^p\, d\mu^2(x,y)
 & & \\
 & & \hskip-48mm \leq \
\frac{p^2}{4}\, \int\!\!\!\int |u(x) - u(y)|^{p-2}\,
\Big(|\nabla u(x)|^2 + |\nabla u(y)|^2\Big)\,d\mu^2(x,y) \\
 & & \hskip-48mm \leq \
\frac{p^2}{4}\, 
\bigg(\int\!\!\!\int |u(x) - u(y)|^p\,d\mu^2(x,y)\bigg)^{\frac{p-2}{p}}
\bigg(\int\!\!\!\int 
\Big(|\nabla u(x)|^2 + |\nabla u(y)|^2\Big)^{\frac{p}{2}}\,
d\mu^2(x,y)\bigg)^{\frac{2}{p}}.
\ene
By Jensen's inequality, the last double integral does not exceed
$$
2^{\frac{p}{2} - 1} \int\!\!\!\int 
\big(|\nabla u(x)|^p + |\nabla u(y)|^p\big)\,d\mu^2(x,y) \,= \,
2^{\frac{p}{2}} \int |\nabla u|^p\,d\mu,
$$
and hence
$$
\lambda_1 
\bigg(\int\!\!\!\int |u(x) - u(y)|^p\,d\mu^2(x,y)\bigg)^{\frac{2}{p}}
 \, \leq \, \frac{p^2}{2}\,
 \bigg(\int |\nabla u|^p\,d\mu\bigg)^{\frac{2}{p}}.
$$
Equivalently,
$$
\int\!\!\!\int |u(x) - u(y)|^p\,d\mu^2(x,y) \, \leq \,
\Big(\frac{p}{\sqrt{2\lambda_1}}\Big)^p\,\int |\nabla u|^p\,d\mu.
$$

If the right integral is finite, so is the left one, thus
$u$ is integrable. Moreover, 
the left integral is greater than or equal to $\int |u(x) - \E\,u(X)|^p\,d\mu(x)$ 
(by Jensen's inequality). 
\qed

\vskip5mm
Let us now connect the Poincar\'e constant with small ball probabilities. 

\vskip5mm
{\bf Corollary 9.3.} {\sl If $\E\,|X|^2 = n$, then
\be
\P\Big\{|X|^2 \leq \frac{1}{4}\,n\Big\} \, \leq \,
3\,e^{- \frac{1}{2}\sqrt{\lambda_1 n}}.
\en
}

\vskip2mm
{\bf Proof.} Note that $\lambda_1 \leq 1$ due to the assumption
$\E\,|X|^2 = n$. Applying (9.2) to the function $u(x) = -|x|$, we have
\be
\P\big\{|X| - \E\, |X| \leq -r\big\} \, \leq \,
3\,e^{-2\sqrt{\lambda_1} r}, \quad r \geq 0.
\en
One can bound $\E\, |X|$ from below by virtue of the Poincar\'e-type 
inequality (9.1) which gives 
$$
n - (\E\,|X|)^2 = \Var(|X|) \leq \, \frac{1}{\lambda_1}.
$$
In the case $\lambda_1 n \geq \frac{4}{3}$, this implies
$\E\,|X| \geq \sqrt{n - \frac{1}{\lambda_1}} \geq \frac{1}{2}\sqrt{n}$. 
Hence, applying (9.6) with $r = \E\, |X| -  \frac{1}{2}\sqrt{n}$, we get
$$
\P\Big\{|X| \leq \frac{1}{2}\sqrt{n}\Big\} \leq  3\,e^{-2 \sqrt{\lambda_1}\,r}.
$$
Here $r \geq \sqrt{n - \frac{1}{\lambda_1}} -  \frac{1}{2}\sqrt{n} \geq 
\frac{1}{4}\sqrt{n}$ under a stronger assumption
$\lambda_1 n \geq \frac{16}{7}$, in which case the above bound yields
the desired inequality (9.5).

It remains to note that (9.5) is fulfilled automatically when 
$\lambda_1 n < \frac{16}{7}$, since then the right-hand
side is greater than 1.
\qed

\vskip5mm
Let us give another version of this statement for convolutions, 
namely, for sums
$$
\Sigma_p = \sum_{k=1}^p\, (X^{(k)} - Y^{(k)}),
$$
where $X^{(k)}$, $Y^{(k)}$ ($1 \leq k \leq p$) are independent copies of 
$X$. One may use the property that the product measure $\mu^{\otimes 2p}$ 
on $(\R^n)^{2p} = \R^{2pn}$ has the same Poincar\'e constant $\lambda_1$ 
as the distribution $\mu$ of $X$. The function
$$
u(x_1,\dots,x_p,y_1,\dots,y_p) \, = \, -
\Big|\sum_{k=1}^p\, (x_k - y_k)\Big|, \qquad x_k,y_k \in \R^n,
$$
has Lipschitz semi-norm $\sqrt{2p}$ with respect to the Euclidean distance 
on $\R^{2pn}$. Therefore, according to Proposition 9.1, it admits 
an exponential inequality
$$
\mu^{\otimes 2p}\{u - m \geq r\} \, \leq \,
3e^{-2\sqrt{\lambda_1}\,r/\sqrt{2p}} \qquad (r > 0),
$$
where $m$ is the $\mu^{\otimes 2p}$-mean of $u$. That is,
\be
\P\big\{|\Sigma_p| - \E\,|\Sigma_p|\leq -r\big\} \, \leq \,
3e^{-2\sqrt{\lambda_1}\,r/\sqrt{2p}}.
\en

By the Poincar\'e-type inequality on the product space, and using 
$\E\,|\Sigma_p|^2 = 2p n$, we have
$$
2pn - (\E\,|\Sigma_p|)^2 \leq \, \frac{2p}{\lambda_1} \leq pn,
$$
where the last inequality holds true when $\lambda_1 n \geq 2$.
In this case, $\E\,|\Sigma_p| \geq \sqrt{pn}$, and applying (9.7) with
$r = (1 - \frac{1}{\sqrt{2}})\sqrt{pn}$, we obtain
$$
\P\Big\{|\Sigma_p| \leq \frac{1}{\sqrt{2}}\sqrt{np}\Big\} \, \leq \,
3\,e^{-(\sqrt{2} - 1) \sqrt{\lambda_1 n}} <
3\,e^{-\frac{1}{3} \sqrt{\lambda_1 n}}.
$$
In the case $\lambda_1 n \leq 2$, this inequality is fulfilled automatically, so,
we arrive at:

\vskip5mm
{\bf Corollary 9.4.} {\sl If $\E\,|X|^2 = n$, then
$$
\P\Big\{|\Sigma_p|^2 \leq \frac{np}{2}\Big\} \, \leq \, 
3\,e^{-\frac{1}{3} \sqrt{\lambda_1 n}}.
$$
}

\vskip2mm
{\bf Remark 9.5.} If the random vector $X$ in $\R^n$ ($n \geq 2$) is isotropic,
then necessarily $\lambda_1 \leq 1$. Indeed, applying the Poincar\'e-type
inequality (9.1) with linear functions $u(x) = \left<x,\theta\right>$, we get
$$
\lambda_1\,\big(1 - \left<a,\theta\right>^2\big) \leq 1, \quad \theta \in \S^{n-1},
$$
where $a = \E X$. Since one may choose $\theta$ to be orthogonal to the vector $a$,
the conclusion follows. The upper bound $\lambda_1 \geq 1$ is also valid 
in dimension $n=1$, as long as $\E X = 0$ (however, we only have
$\lambda_1 \leq 1/\Var(X)$ without the mean zero assumption).


\vskip7mm
\section{{\bf The Case of Non-symmetric Distributions}}
\setcounter{equation}{0}

\vskip2mm
\noindent
In order to extend the bound
\be
\E_\theta\, \rho(F_\theta,\Phi) \leq \frac{c\,\log n}{n}\, \Lambda
\en
to the case where the distribution of $X$ is not necessarily symmetric 
about the origin, we need to employ more sophisticated results reflecting 
the size of the linear part of the characteristic functions $f_\theta(t)$ 
in $L^2(\s_{n-1})$ with respect to the variable $\theta$. This may be 
achieved at the expense of a certain term that has to be added to the 
right-hand side in (10.1). More precisely, we derive the following:

\vskip5mm
{\bf Proposition 10.1.} {\sl Given an isotropic random vector
$X = (X_1,\dots,X_n)$ in $\R^n$,
\be
c\,\E_\theta\, \rho(F_\theta,\Phi) \, \leq \, 
\frac{\log n}{n}\, \Lambda + \Big(\frac{\log n}{n}\Big)^{1/4}
\bigg(\E\, \frac{\left<X,Y\right>}{\sqrt{|X|^2 + |Y|^2}}\bigg)^{1/2},
\en
where $Y$ is an independent copy of $X$. 
}

\vskip5mm 
The ratio $\left<X,Y\right>/\sqrt{|X|^2 + |Y|^2}$ is understood to be 
zero in the case $X = Y = 0$. Note that the last expectation in (10.2)
is non-negative which follows from the representation
$$
\E\, \frac{\left<X,Y\right>}{\sqrt{|X|^2 + |Y|^2}} \, = \,
\frac{2}{\sqrt{\pi}} \int_0^\infty \sum_{k = 1}^n 
\Big(\E X_k \, e^{-|X|^2 r^2}\Big)^2\,dr.
$$
If the distribution of $X$ is symmetric, this expectation is vanishing, 
and we return to (10.1).

Returning to Proposition 6.1, define the random variables
$$
R^2 = \frac{|X|^2 + |Y|^2}{2n} \ \ (R \geq 0), \qquad 
U = \frac{|X|^2}{n}, \ V = \frac{|Y|^2}{n},
$$
and recall that the squared $L^2$-norm of the linear part of the
characteristic function $f_\theta(t)$ of the weighted sums
$\left<X,\theta\right>$ admits an asymptotic representation
\be
I(t) = \frac{t^2}{n}\, \E \left<X,Y\right>
\Big(1 - \frac{(U^2 + V^2)\, t^4 - 8R^2 t^2}{4n}\Big)\,
e^{-R^2 t^2} + O(t^2 n^{-5/2}).
\en

\vskip5mm
{\bf Lemma 10.2.} {\sl If $X$ is isotropic, then,
putting $T_0 = 4\sqrt{\log n}$, we have
\be
\int_0^{T_0} \frac{I(t)}{t^2} \,dt \, \leq \,
\frac{c}{n}\,\E\,\frac{\left<X,Y\right>}{R} +
O\big(\Lambda^2 n^{-2}\big).
\en
}

\vskip2mm
{\bf Proof.} Introduce the events 
$A = \{R \leq \frac{1}{2}\}$ and $B = \{R > \frac{1}{2}\}$. 
From (10.3),
\bee
\int_0^{T_0} \frac{I(t)}{t^2}\,dt
 & = &
\frac{1}{n}\,\E\, 
\left<X,Y\right> \int_0^{T_0} e^{-R^2 t^2}\,dt
 +
\frac{2}{n^2}\,\E\, 
\left<X,Y\right> \int_0^{T_0} R^2 t^2 e^{-R^2 t^2}\,dt \\
 & & - \ 
\frac{1}{4n^2}\,\E\,\left<X,Y\right> \int_0^{T_0} 
(U^2 + V^2)\, t^4 e^{-R^2 t^2}\,dt + O(n^{-2}).
\ene
After the change of the variable $Rt = s$ (assuming without loss 
of generality that $R>0$) and putting $T_1 = RT_0$, the above 
is simplified to
\bee
\int_0^{T_0} \frac{I(t)}{t^2}\,dt
 & = &
\frac{1}{n}\,\E\,\frac{\left<X,Y\right>}{R} 
\int_0^{T_1} e^{-s^2}\,ds +
\frac{2}{n^2}\,\E\,\frac{\left<X,Y\right>}{R} 
\int_0^{T_1} s^2 e^{-s^2}\,ds \\
 & & - \ 
\frac{1}{4n^2}\,
\E\,\frac{\left<X,Y\right>}{R}\,\frac{U^2 + V^2}{R^4}
\int_0^{T_1} s^4 e^{-s^2}\,ds + O(n^{-2}).
\ene 

At the expense of a small error, integration here may be 
extended from the interval $[0,T_1]$ to the whole half-axis 
$(0,\infty)$. To see this, one can use the estimates
$$
\int_{T_1}^\infty  e^{-s^2}\,ds < \int_{T_1}^\infty  s^2\,e^{-s^2}\,ds <
\int_{T_1}^\infty  s^4 e^{-s^2}\,ds < c\,e^{-T_1^2/2} \qquad
(T_1 > 1), 
$$
together with
\be
\Big|\frac{\left<X,Y\right>}{R}\Big| \leq 
\frac{|X| \, |Y|}{R} \leq \frac{|X|^2 + |Y|^2}{2R} = Rn.
\en
As was already noted in (7.4),
\begin{eqnarray}
\P(A) 
 & = & 
\P\Big\{|X|^2 + |Y|^2 \leq \frac{n}{2}\Big\} \nonumber \\
 & \leq & 
\P\Big\{|X|^2 \leq \frac{n}{2}\Big\} \,
\P\Big\{|Y|^2 \leq \frac{n}{2}\Big\} \, \leq \, 
\frac{16\Lambda^2}{n^2}.
\end{eqnarray}
Since on the set $B$, we have  
$T_1^2 = 16 R^2 \log n > 4\,\log n$, and due to $\E R^2 = 1$, 
it follows that
\bee
\E\,R\, e^{-T_1^2/2} 
 & = &
\E\,R\, e^{-T_1^2/2}\,1_A + \E\,R\, e^{-T_1^2/2}\,1_B \\
 & \leq &
\frac{1}{2}\,\P(A) + \frac{1}{n^2}\,\E R
 \, \leq \,
\frac{c\Lambda^2}{n^2},
\ene
where we used the lower bound $\Lambda \geq \frac{1}{2}$. Hence
$$
\E\,\frac{|\left<X,Y\right>|}{R} \int_{T_1}^\infty  e^{-s^2}\,ds \, \leq \,
n\,\E R e^{-T_1^2/2} \, \leq \, \frac{c\Lambda^2}{n}.
$$
By a similar argument,
$$
\E\,\frac{|\left<X,Y\right>|}{R} 
\int_{T_1}^\infty s^2 e^{-s^2}\,ds \, \leq \,
cn\,\E R\,e^{-T_1^2/2} \, \leq \, \frac{c\Lambda^2}{n}.
$$
Using 
$$
\frac{U^2 + V^2}{R^4} \, = \, \frac{4\,(U^2 + V^2)}{(U+V)^2} \, \leq \, 4,
$$
we also have
$$
\E\,\frac{|\left<X,Y\right>|}{R}\,\frac{U^2 + V^2}{R^4} 
\int_{T_1}^\infty s^4 e^{-s^2}\,ds \, \leq \, 
cn\, \E R\,e^{-T_1^2/2} \, \leq \, \frac{c\Lambda^2}{n}.
$$

Thus, extending the integration to the positive half-axis, we get
\bee
\int_0^{T_0} \frac{I(t)}{t^2}\,dt 
 & = &
\frac{c_1}{n}\,\E\,\frac{\left<X,Y\right>}{R} +
\frac{c_2}{n^2}\,\E\,\frac{\left<X,Y\right>}{R} \\
 & & - \
\frac{c_3}{n^2}\,\E\,
\frac{\left<X,Y\right>}{R}\, \frac{U^2 + V^2}{R^4}
 + O\big(\Lambda^2 n^{-2}\big)
\ene
with some absolute constants $c_j>0$. Moreover, using the identity
$$
\frac{U^2 + V^2}{R^4} = 2 + \frac{(U - V)^2}{2R^4} = 
2 + 2\,\frac{(U - V)^2}{(U+V)^2}
$$
and recalling that  $\E\,\frac{\left<X,Y\right>}{R} \geq 0$, it follows that, 
with some other positive absolute constants
\be
\int_0^{T_0} \frac{I(t)}{t^2}\,dt \, \leq \,
\frac{c_1}{n}\,\E\,\frac{\left<X,Y\right>}{R} -
\frac{c_2}{n^2}\,\E\,\frac{\left<X,Y\right>}{R}\,
\frac{(U - V)^2}{(U + V)^2} + O\big(\Lambda^2 n^{-2}\big).
\en

To get rid of the last expectation (by showing that it is bounded
by a dimension free quantity), first note that, by (10.5),
the expression under this expectation is bounded in absolute value by
$Rn$. Hence, applying Cauchy's inequality and using $\E R^2 = 1$, from
(10.6) we obtain that
\begin{eqnarray}
\E\,\Big|\frac{\left<X,Y\right>}{R}\Big|\, \frac{(U-V)^2}{(U+V)^2}\,1_A
 & \leq & 
\E\,\Big|\frac{\left<X,Y\right>}{R}\Big|\,1_A \nonumber \\
 &  \leq &
n\,\E R\,1_A  \, \leq \, n\sqrt{\P(A)} \, \leq \, 4\Lambda.
\end{eqnarray}
Turning to the complementary set, note that on $B$, we have 
$|\frac{\left<X,Y\right>}{R}| \leq 2\,|\left<X,Y\right>|$, while
$$
\frac{(U - V)^2}{(U+V)^2} \, \leq \, \frac{|U - V|}{U+V} \, = \,
\frac{|U - V|}{2R^2} \, \leq \, 2\,|U - V|.
$$ 
Hence, by Cauchy's inequality, and using $\E \left<X,Y\right>^2 = n$,
\bee
\E\,\Big|\frac{\left<X,Y\right>}{R}\Big|\, \frac{(U - V)^2}{(U+V)^2}\,1_B
 & \leq & 
4\,\E\,|\left<X,Y\right>|\,|U - V| \\
 & \leq &
4\sqrt{n}\,\sqrt{\E\,(U-V)^2} \, = \, 4\sqrt{2}\,\sigma_4 \, \leq \,
4\sqrt{2\Lambda}.
\ene
Combining this bound with (10.8), we finally obtain that
$$
\E\,\Big|\frac{\left<X,Y\right>}{R}\Big|\, \frac{(U-V)^2}{(U+V)^2} \, \leq \, 
c\Lambda.
$$
As a result, we arrive in (10.7) at the bound (10.4).

\vskip5mm
{\bf Proof of Proposition 10.1.}  We employ the bound (7.7) of 
Lemma 7.2 which was stated with $T_0 = 4\sqrt{\log n}$. Using Cauchy's 
inequality and applying (10.4), it gives
\bee
c\,\E_\theta\, \rho(F_\theta,\Phi) 
 & \leq &
\frac{\log n}{n}\, \Lambda + \int_0^{T_0} \frac{\sqrt{I(t)}}{t}\,dt \\
 & \leq &
\frac{\log n}{n}\, \Lambda +
\sqrt{T_0} \left(\int_0^{T_0} \frac{I(t)}{t^2}\,dt\right)^{1/2} \\
 & \leq &
\frac{\log n}{n}\, \Lambda +
c' \sqrt{T_0}\, \bigg(\frac{1}{n}\,\E\,\frac{\left<X,Y\right>}{R} +
\frac{\Lambda^2}{n^2}\bigg)^{1/2}.
\ene
Simplifying the expression on the right-hand side, we arrive at (10.2).
\qed

\vskip7mm
\section{{\bf The estimate on average}}
\setcounter{equation}{0}

\vskip2mm
\noindent
Let us rewrite the bound (10.2) as
\be
c\,\E_\theta\, \rho(F_\theta,\Phi) \, \leq \, 
\frac{\log n}{n}\, \Lambda + \frac{(\log n)^{1/4}}{\sqrt{n}}
\Big(\E\, \frac{\left<X,Y\right>}{R}\Big)^{1/2},
\en
where $R^2 = \frac{1}{2n}\,(|X|^2 + |Y|^2)$, $R \geq 0$, and where 
$Y$ is an independent copy of $X$. In the next step, we are going 
to simplify the last expectation in terms of $\lambda_1$.
Note that, under our standard assumptions as in Proposition 10.1, 
$$
\E R^2 = 1, \quad \Var(R^2) = \frac{\sigma_4^2}{2n} \leq \frac{\Lambda}{2n}. 
$$
Hence, with high probability the ratio $\frac{\left<X,Y\right>}{R}$ 
is almost $\left<X,Y\right>$ which in turn has zero expectation, 
as long as $X$ has mean zero. However, in general it is not clear 
whether or not this approximation is sufficient to make further 
simplification. Nevertheless, the approximation $R^2 \sim 1$ is indeed 
sufficiently strong, for example, in the case where the distribution 
$\mu$ of $X$ satisfies the Poincar\'e-type inequality (1.3).

\vskip5mm
{\bf Lemma 11.1.} {\sl Let $X$ be an isotropic random vector in $\R^n$ 
with mean zero and a positive Poincar\'e constant $\lambda_1$, and let
$Y$ be an independent copy of $X$. Then
\be
\E\, \frac{\left<X,Y\right>}{R} \, \leq \, \frac{c}{\lambda_1^2\, n}.
\en
}

\vskip2mm
Applying (11.2) in (11.1) and using $\Lambda \leq 4/\lambda_1$
(cf. \cite{B-C-G5}, Proposition 3.4), we get an estimate on average
$$
c\,\E_\theta\, \rho(F_\theta,\Phi) \, \leq \, 
\frac{\log n}{n}\, \frac{1}{\lambda_1} + 
\frac{(\log n)^{1/4}}{\sqrt{n}} \, \frac{1}{\lambda_1\sqrt{n}},
$$
thus proving the relation (1.6).

\vskip5mm
{\bf Proof of Lemma 11.1.} Without loss of generality, assume that $R > 0$ a.s.
Put $\delta_n = \frac{1}{\lambda_1 n}$.

We apply the Poincar\'e-type inequality for 
the product measure $\mu \otimes \mu$,
\be
\int\!\!\!\int |u(x,y)|^2\,d\mu(x)\,d\mu(y) \, \leq \, 
\frac{1}{\lambda_1} 
\int\!\!\!\int |\nabla u(x,y)|^2\,d\mu(x)\,d\mu(y),
\en
which holds true for any smooth function $u$ on $\R^n \times \R^n$
with $(\mu \otimes \mu)$-mean zero. Moreover,
according to the inequality (9.3), for any $p \geq 2$,
\be
\int\!\!\!\int |u(x,y)|^p\,d\mu(x)\,d\mu(y) \, \leq \, 
\frac{p^p}{(2\lambda_1)^{p/2}} 
\int\!\!\!\int |\nabla u(x,y)|^p\,d\mu(x)\,d\mu(y).
\en
By Corollary 9.3 applied in $\R^{2n}$ to the random vector $(X,Y)$, 
it also follows that the event $A = \{R \leq \frac{1}{2}\}$ has probability
$$
\P(A) \leq 3e^{-\sqrt{\lambda_1 n/2}}.
$$

Using
\be
|\left<X,Y\right>| \leq R^2 n,
\en 
cf. (10.5), we have
\be
\E\, \frac{|\left<X,Y\right>|}{R}\,1_A \, \leq \, n\,\E R\,1_A \, \leq \,
\frac{n}{2}\ \P(A) \, \leq \,
\frac{3n}{2}\, e^{-\sqrt{\lambda_1 n/2}} \, \leq \, 
\frac{c}{\lambda_1^2\, n}.
\en
Similarly,
$$
\E\,|\left<X,Y\right>|\,1_A \, \leq \, \frac{n}{4}\ \P(A)
 \, \leq \, \frac{c}{\lambda_1^2\, n},
$$
and since $X$ has mean zero, for the complementary set 
$B = \{R > \frac{1}{2}\}$ we have the same bound
$$
\big|\,\E\left<X,Y\right> 1_B\big| \leq 
\frac{c}{\lambda_1^2\, n}.
$$
Using once more (11.5), on the set $A$ we also have
$$
\E\,|\left<X,Y\right>
|\,R^2\,1_A \leq \frac{n}{16}\,\P(A)\leq \frac{c}{\lambda_1^2 n}
$$
and
$$
\E\,|\left<X,Y\right>|\,R^4\,1_A \leq \frac{n}{64}\,\P(A)\leq 
\frac{c}{\lambda_1^2\, n}.
$$

Now, consider the function $w(\ep) = (1 + \ep)^{-1/2}$ 
on the half-axis $\ep \geq -\frac{3}{4}$. By Taylor's formula, 
for some point $\ep_1$ between $-\frac{3}{4}$ and $\ep$,
$$
w(\ep) \, = \, 
1 - \frac{1}{2}\,\ep + \frac{3}{8}\,\ep^2 -
\frac{5}{16}\,(1 + \ep_1)^{-7/2}\,\ep^3 \, = \, 
1 - \frac{1}{2}\,\ep + \frac{3}{8}\,\ep^2 - \beta \ep^3
$$
with some $0 \leq \beta \leq 40$. Putting 
$\ep = R^2 - 1$, we then get on the set $B$
\bee
\frac{\left<X,Y\right>}{R} 
 & = & 
\left<X,Y\right> - \frac{1}{2} \left<X,Y\right>(R^2 - 1) + 
\frac{3}{8}\left<X,Y\right>(R^2 - 1)^2 - 
\beta \left<X,Y\right> (R^2 - 1)^3 \\ 
 & = & 
\frac{15}{8}\left<X,Y\right> - \frac{5}{4}\left<X,Y\right> R^2 + 
\frac{3}{8}\left<X,Y\right> R^4 - 
\beta \left<X,Y\right> (R^2 - 1)^3.
\ene
By the independence of $X$ and $Y$, and due to the mean zero assumption, 
$\E \left<X,Y\right> = \E \left<X,Y\right> R^2 = 0$.
Hence, writing $1_B = 1 - 1_A$, we have
\bee
\E\, \frac{\left<X,Y\right>}{R}\,1_B
 & = & 
-\frac{15}{8}\,\E \left<X,Y\right> 1_A + 
\frac{5}{4}\,\E \left<X,Y\right> R^2\,1_A -
\frac{3}{8}\,\E \left<X,Y\right> R^4\,1_A \\
 & & + \  
\frac{3}{8}\,\E \left<X,Y\right> R^4 -
\beta\,\E \left<X,Y\right>(R^2 - 1)^3\,1_B.
\ene
Here, the first three expectations on the right-hand side 
do not exceed in absolute value a multiple of $\frac{1}{\lambda_1^2 n}$. 
Hence, using the previous bound (11.6), we get
\be
\E\, \frac{\left<X,Y\right>}{R} \, = \,
\frac{c_1}{\lambda_1^2\, n} +  
\frac{3}{8}\,\E \left<X,Y\right> R^4 +
c_2\,\E\, |\left<X,Y\right>|\,|R^2 - 1|^3,
\en
where the quantities $c_1$ and $c_2$ are bounded by 
an absolute constant.

By Cauchy's inequality, the square of the last expectation does not exceed, 
$$
\E \left<X,Y\right>^2\, \E\,(R^2 - 1)^6 \, = \, 
n \, \E\,(R^2 - 1)^6.
$$
In turn, the latter expectation may be bounded by virtue of 
the inequality (11.4) applied with $p=6$ to the function 
$u(x,y) = \frac{1}{2n}\,(|x|^2 + |y|^2) - 1$. Since
$$
|\nabla u(x,y)|^2 = |\nabla_x u(x,y)|^2 + |\nabla_y u(x,y)|^2 =
\frac{|x|^2 + |y|^2}{n^2},
$$
it gives
\be
\E\,(R^2 - 1)^6 \, \leq \, \frac{c}{\lambda_1^3\, n^3}\,\E R^6.
\en
On the other hand, the Poincar\'e-type inequality easily yields
the bound $\E\,R^6 \leq c/\lambda_1^3$. 
However, in this step a more accurate estimation is required. Write
$$
R^6 = (R^2 - 1)^3 + 3\,(R^2 - 1)^2 + 3\,(R^2 - 1) + 1,
$$
so that
\be
\E R^6 = \E\,(R^2 - 1)^3 + 3\,\E\,(R^2 - 1)^2 + 1.
\en
By (11.3) with the same function $u$, we have 
$$
\E\,(R^2 - 1)^2 \leq \frac{2}{\lambda_1 n}\,\E R^2 = 2\delta_n,
$$ 
while (11.4) with $p=3$ gives 
$$
\E\,|R^2 - 1|^3 \, \leq \, 27\,\delta_n^{3/2}\,\E\,|R|^3.
$$
Putting $x^2 = \E R^6$ ($x > 0$) and using
$\E\,|R|^3 \leq x$, we therefore get from (11.9) that
$$
x^2 \, \leq \, 27\,\delta_n^{3/2} x + 6\delta_n + 1.
$$
This quadratic inequality is easily solved to yield
$x^2 \leq c\,(\delta_n + 1)^3$.
One can now apply this bound in (11.8) to conclude that
$$
\E\,(R^2 - 1)^6 \, \leq \, 
\frac{c}{\lambda_1^3\, n^3}\,(\delta_n + 1)^3.
$$

This implies
$$
\E \left<X,Y\right>^2\, \E\,(R^2 - 1)^6 \, \leq \, 
\frac{c}{\lambda_1^3\, n^2}\,(\delta_n + 1)^3,
$$
which allows us to simplify the representation (11.7) 
to the form
\be
\E\, \frac{\left<X,Y\right>}{R} \, = \,
\frac{c_1}{\lambda_1^2\, n} +  
\frac{c_2}{\lambda_1^{3/2}\, n}\,(\delta_n + 1)^{3/2} +
\frac{3}{8}\,\E \left<X,Y\right> R^4,
\en
where the new quantity $c_2$ is bounded by an absolute
constant.

We are left with the estimation of $\E \left<X,Y\right> R^4$. Since 
$\E \left<X,Y\right> |X|^4 = \E \left<X,Y\right> |Y|^4 = 0$,
it follows that
$$
\E \left<X,Y\right> R^4 \, = \, 
\frac{1}{2n^2}\,\E \left<X,Y\right> |X|^2\, |Y|^2
 \, = \, 
\frac{1}{2n^2}\,\big|\,\E\, |X|^2 X\,\big|^2.
$$
Here the latter expectation is understood in the usual vector 
sense. That is, in terms of the components in $X = (X_1,\dots,X_n)$
defined on a probability space $(\Omega,\P)$, we have
$$
\E\, |X|^2 X = (a_1,\dots,a_n), \qquad a_k = \E\,|X|^2 X_k = 
\E\,(|X|^2 - n)\, X_k.
$$
Since the collection $\{X_1,\dots,X_n\}$ appears as an orthonormal system 
in the Hilbert space $L^2(\Omega,\P)$, the numbers $a_k$ represent the (Fourier) 
coefficients for the projection of the random variable
$|X|^2 - n$ onto the span of $X_k$'s. Hence, by Bessel's inequality,
$$
\big|\,\E\, |X|^2 X\,\big|^2 \, = \, \sum_{k=1}^n a_k^2 \, \leq \,
\big\|\,|X|^2 - n\big\|_{L^2(\Omega,\P)}^2 \, = \,
\Var(|X|^2) \, = \, n\,\sigma_4^2(X) \, \leq \, \frac{4n}{\lambda_1},
$$
so that
$$
\E \left<X,Y\right> R^4 \,\leq \, \frac{2}{\lambda_1 n}
$$
In view of the upper bound $\lambda_1 \leq 1$ (Remark 9.5), 
the expectation in (11.10) is thus dominated by the first term, and we arrive at
$$
\E\, \frac{\left<X,Y\right>}{R} \, \leq \, 
\frac{c}{\lambda_1^2\, n} +  
\frac{c}{\lambda_1^{3/2}\, n}\,\Big(\frac{1}{\lambda_1 n} + 1\Big)^{3/2}.
$$

If $\lambda_1 \geq n^{-3/2}$, the first term on the right-hand side 
dominates the second one, and we arrive at the desired inequality (11.2). 
In the other case, we have $\frac{1}{\lambda_1^2\, n} \geq n^2$, and 
then (11.2) holds true as well, by (11.5), since $\E R \leq 1$.
\qed

\vskip7mm
\section{{\bf Proof of Theorem 1.1}}
\setcounter{equation}{0}

\vskip2mm
\noindent
Let us now derive the stronger inequality (1.7). 
With parameters $T_0 = 4\sqrt{\log n}$ and $T = T_0 n$,
the bound (7.2) of Lemma 7.1 is simplified to
\be
c\,\rho(F_\theta,\Phi) \, \leq \,
\int_0^{T_0} \frac{|f_\theta(t) - f(t)|}{t}\,dt +
L(\theta)  + \frac{\log n}{n}\,\Lambda,
\en
where 
$
L(\theta) = \int_{T_0}^T \frac{|f_\theta(t)|}{t}\,dt.
$
Combining Corollary 9.4 with Lemma 8.2, we obtain that
$$
\E_\theta\,L(\theta)^{2p} \, \leq \, (c\log n)^{2p}\, 
\Big(p^{2p}\, n^{-2p} + 
e^{-\frac{1}{3}\sqrt{\lambda_1 n}}\,\Big)
$$
for any integer $p \geq 1$.
One can simplify this bound, by using the inequality 
$e^{-x} \leq (\frac{4p}{ex})^{4p}$ ($x>0$). Since
$\lambda_1 \leq 1$ (as was explained above), it follows that
$$
\big(\E_\theta\,L(\theta)^{2p}\big)^{1/2p} \, \leq \, 
\frac{c\log n}{n}\,\lambda_1^{-1}\, p^2.
$$
This inequality is readily extended to all real $p \geq 1/2$.
Replacing here $2p$ with $p$ we get a similar bound
$$
\big(\E_\theta\,L(\theta)^p\big)^{1/p} \, \leq \, 
\frac{c\log n}{n}\,\lambda_1^{-1}\, p^2,
$$
which holds for all real $p \geq 1$.
Now, by Markov's inequality,
$$
\s_{n-1}\Big\{L(\theta) \geq 
\frac{ce\log n}{n}\,\lambda_1^{-1} r\Big\} \leq 
\frac{p^{2p}}{(er)^p}, \qquad r \geq 1.
$$
Choosing $p = \sqrt{r}$, we thus have
\be
\s_{n-1}\Big\{L(\theta) \geq 
\frac{ce\log n}{n}\,\lambda_1^{-1} r\Big\} \, \leq \, e^{-\sqrt{r}}.
\en

It is time to involve Lemma 8.1. 
First, from Lemmas 10.2 and 11.1, it follows that
\bee
\int_0^{T_0} \frac{\sqrt{I(t)}}{t}\,dt
 & \leq &
\sqrt{T_0} \left(\int_0^{T_0} \frac{I(t)}{t^2}\,dt\right)^{1/2} \\
 & \leq &
c \sqrt{T_0}\, \bigg(\frac{1}{n}\,\E\,\frac{\left<X,Y\right>}{R} +
\frac{\Lambda^2}{n^2}\bigg)^{1/2} \, \leq \
\frac{c'}{\lambda_1 n}\,(\log n)^{1/4},
\ene
where on the last step we used $\Lambda \leq \frac{4}{\lambda_1}$. 
Hence, by Lemma 8.1,
$$
\s_{n-1}\bigg\{\int_0^{T_0} \frac{|f_\theta(t) - f(t)|}{t}\,dt 
\geq \frac{c\log n}{\lambda_1 n}\,r\bigg\} \, \leq \, 2\,e^{-r}.
$$
Being combined with (12.2) and applied in (12.1), this bound
leads to the desired inequality
\be
\s_{n-1}\Big\{\rho(F_\theta,\Phi) \geq 
\frac{c\log n}{n}\,\lambda_1^{-1} r\Big\} \, \leq \, 3\,e^{-\sqrt{r}},
\en
which also holds for $r < 1$ (when the right-hand side is greater
than 1). Here, the constant 3 may be replaced with 2 by rescaling
the variable $r$, and then we arrive at (1.7).
\qed

\vskip5mm
{\bf Corollary 12.1.} {\sl Let $X$ be an isotropic random vector in $\R^n$ 
with mean zero and a positive Poincar\'e constant $\lambda_1$.
For any $\beta>0$, with $\s_{n-1}$-probability at most $3n^{-\beta}$
we have
$$
\rho(F_\theta,\Phi) \leq \frac{c\beta^2\,(\log n)^3}{n}\,\lambda_1^{-1}.
$$ 
}

\vskip2mm
Indeed, although the estimate (1.7) implies the bound on average (1.6),
it is only effective for $r \geq (\log n)^2$. For the values 
$r = (\beta \log n)^2$, (12.3) provides a polynomial bound
$$
\s_{n-1}\Big\{\rho(F_\theta,\Phi) \geq 
\frac{c\beta^2\,(\log n)^3}{n}\,\lambda_1^{-1}\Big\} \, \leq \, 
3n^{-\beta}.
$$
In other words, for a sufficiently large number $A$, with high 
$\s_{n-1}$-probability
$$
\rho(F_\theta,\Phi) \leq \frac{A\,(\log n)^3}{n}\,\lambda_1^{-1}.
$$


\begin{thebibliography}{B}
\itemsep=6pt
\small

\vskip3mm
\bibitem{A-S}
S. Aida, and D. Stroock. Moment estimates derived from Poincar\'e 
      and logarithmic Sobolev inequalities. Math. Res. Lett. 1 (1994), 
			no. 1,  75--86. 

\bibitem{A-B-P}
M. Anttila, K. Ball, and I. Perissinaki. The central limit problem for 
      convex bodies. Trans. Amer. Math. Soc. 355 (2003), no. 12, 4723--4735.

\bibitem{Bat}
H. Bateman. Higher transcendental functions, Vol. II. McGraw-Hill 
      Book Company, Inc., 1953, 396 pp.

\bibitem{B1}
S. G. Bobkov. Remarks on the Gromov-Milman inequality. (Russian) 
      Vestn. Syktyvkar. Univ. Ser. 1 Mat. Mekh. Inform. No. 3 (1999), 15--22.

\bibitem{B2}
S. G. Bobkov. On concentration of distributions of random weighted sums. 
      Ann. Probab. 31 (2003), 195--215.

\bibitem{B3} 
S. G. Bobkov. On a theorem of V. N. Sudakov on typical distributions. 
      Zap. Nauchn. Sem. S.-Peterburg. Otdel. Mat. Inst. Steklov. (POMI) 
			368 (2009), 59--74, 283.

\bibitem{B4}
S. G. Bobkov. Closeness of probability distributions in terms of 
      Fourier-Stieltjes transforms. Russian Math. Surveys, 
			vol. 71, issue 6, (2016), 1021--1079. Translated from: 
			Uspekhi Matem. Nauk, vol. 71, issue 6 (432), (2016), 37--98.

\bibitem{B5}
S. G. Bobkov. Edgeworth corrections in randomized central limit theorems.
      Geometric Aspects of Functional Analysis, 2256 (2020), 71--97.			

\bibitem{B-C-G1}
S. G. Bobkov, G. P. Chistyakov, and F. G\"otze. Stability problems in 
      Cram\'er-type characterization in case of i.i.d. summands. 
			Theory Probab. Appl. 57 (2013), no. 4, 568--588. 

\bibitem{B-C-G2}
S. G. Bobkov, G. P. Chistyakov, and F. G\"otze. Second-order concentration on 
      the sphere. Commun. Contemp. Math. 19 (2017), no. 5, 1650058, 20 pp.

\bibitem{B-C-G3}
S. G. Bobkov, G. P. Chistyakov, and F. G\"otze. Gaussian mixtures and normal 
      approximation for V. N. Sudakov's typical distributions. 
      Zap. Nauchn. Sem. S.-Peterburg. Otdel. Mat. Inst. Steklov. 
			(POMI) 457 (2017), Veroyatnost i Statistika. 25, 37--52; 
			reprinted in J. Math. Sci. (N.Y.) 238 (2019), no. 4, 366--376.

\bibitem{B-C-G4}
S. G. Bobkov, G. P. Chistyakov, and F. G\"otze. Berry-Esseen bounds for typical
      weighted sums. J. Electron. Probab.  23 (2018), no. 92, 1--22.

\bibitem{B-C-G5}
S. G. Bobkov, G. P. Chistyakov, and F. G\"otze. Normal approximation for 
      weighted sums under a second order correlation condition.
			Ann. Probab. 48 (2020), no. 3, 1202--1219. 

\bibitem{B-U}
A. A. Borovkov, and S. A. Utev. On an inequality and a characterization of 
      the normal distribution connected with it. Probab. Theory 
      Appl. 28 (1983), 209--218.

\bibitem{B-G-V-V}
S. Brazitikos, A. Giannopolous, P. Valettas, and B.-H. Vritsious. Geometry of
      isotropic convex bodies. Mathematical Surveys and Monographs 196. 
			Amer. Math. Soc., Providence, 2014.		

\bibitem{C}
G. P. Chistyakov. A remark on a theorem of N. A. Sapogov on the stability 
      of decompositions of a normal distribution. (Russian) Operators in function 
      spaces and problems in function theory (Russian), 108--116, 147, 
      ``Naukova Dumka'', Kiev, 1987. 

\bibitem{E}
R. Eldan. Thin shell implies spectral gap up to polylog via a stochastic 
      localization scheme. Geom. Funct. Anal. 23  (2013), no. 2, 532--569. 

\bibitem{E-K}
R. Eldan, and B. Klartag. Pointwise estimates for marginals of convex bodies. 
      J. Funct. Anal. 254 (2008), 2275--2293.

\bibitem{G-S}
L. Goldstein, and Q.-M. Shao. Berry-Esseen bounds for projections of 
      coordinate symmetric random vectors. Electron. Commun. Probab. 
			14 (2009), 474--485. 

\bibitem{G-M}
G. Gromov, and V. D. Milman. A topological application of the isoperimetric
      inequality. Amer. J. Math. 105 (1983), 843--854.

\bibitem{G}
U. Grupel. Remarks on the central limit theorem for non-convex bodies. 
      Geometric aspects of functional analysis, 183--198, 
			Lecture Notes in Math., 2116, Springer, Cham, 2014. 

\bibitem{J-L-V}
H. Jiang, Y. T. Lee, and S. S. Vempala. A generalized central limit conjecture 
      for convex bodies. To appear in: Geometric aspects of functional analysis,
			2020. Also: arXiv:1909.13127.
			
\bibitem{K1}
B. Klartag. A central limit theorem for convex sets. Invent. Math. 
      168 (2007), no. 1, 91--131.			
			
\bibitem{K2}
B. Klartag. Power-law estimates for the central limit theorem for convex sets. 
      J. Funct. Anal. 245 (2007), 284--310.	

\bibitem{K3}					
B. Klartag. A Berry–Esseen type inequality for convex bodies with 
      an unconditional basis. Probab. Theory Related Fields 145 (2009), 1--33.					

\bibitem{K-S}
B. Klartag, and S. Sodin. Variations on the Berry-Esseen theorem. 
      Teor. Veroyatn. Primen. 56 (2011), no. 3, 514--533;  
  		reprinted in: Theory Probab. Appl. 56 (2012), no. 3, 403--419.

\bibitem{L1}
M. Ledoux. On Talagrand's deviation inequalities for product measures. 
      ESAIM Probab. Statist. 1 (1995/97), 63--87. 

\bibitem{L2}
M. Ledoux. Concentration of measure and logarithmic Sobolev inequalities.
      S\'eminaire de Probabilit\'es XXXIII. Lect. Notes in Math.
      1709 (1999), 120--216, Springer.

\bibitem{M-M}
M. W. Meckes, and E. Meckes. The central limit problem for random vectors 
      with symmetries. J. Theoret. Probab. 20 (2007), no. 4, 697--720.

\bibitem{M}
M. W. Meckes. Gaussian marginals of convex bodies with symmetries. 
      Beitr\"age Algebra Geom. 50 (2009), no. 1, 101--118. 

\bibitem{P1} 
V. V. Petrov. Sums of independent random variables. Translated from the Russian 
      by A. A. Brown. Ergebnisse der Mathematik und ihrer Grenzgebiete, Band 82. 
      Springer--Verlag, New York--Heidelberg, 1975. x+346 pp.

\bibitem{P2} 
V. V. Petrov. Limit theorems for sums of independent random variables (Russian), 
      Nauka, Moscow, 1987. 318 pp.   

\bibitem{Sa}
N. A. Sapogov. The problem of stability for a theorem of Cram\'er. (Russian) 
      Vestnik Leningrad. Univ. 10 (1955), no. 11, 61--64. 

\bibitem{So} 
S. Sodin. Tail-sensitive Gaussian asymptotics for marginals of concentrated 
      measures in high dimension. Geometric aspects of functional analysis, 
			271--295, Lecture Notes in Math., 1910, Springer, Berlin, 2007. 

\bibitem{Su}   
V. N. Sudakov. Typical distributions of linear functionals in finite-dimensional 
      spaces of high dimension. (Russian) Soviet Math. Dokl. 19 (1978), 1578--1582; 
	    translation in: Dokl. Akad. Nauk SSSR, 243 (1978), no. 6, 1402--1405.  


			       
\end{thebibliography}
\end{document}